\global\def\extension{1}      %    1 = extended %    0 = short

\catcode`\@=11

\newif\if@fewtab\@fewtabtrue
\def\ps@draft{\let\@mkboth\@gobbletwo \def\@oddhead{}
    \def\@oddfoot{\hbox to 7 cm{\tiny \versionno
       \hfil}\hskip -7cm\hfil\rm\thepage \hfil}
    \def\@evenhead{}\let\@evenfoot\@oddfoot}

\def\draftcite#1{\ifnum\draftcontrol=1#1\else{}\fi}
\def\@lbibitem[#1]#2{\item{}\hskip -3cm \hbox to 2cm
{\hfil$\scriptstyle\draftcite{#2}$}\hskip
1cm[\@biblabel{#1}]\if@filesw
     {\def\protect##1{\string ##1\space}\immediate
      \write\@auxout{\string\bibcite{#2}{#1}}}\fi\ignorespaces}

\def\@bibitem#1{\item\hskip -3cm \hbox to 2cm
{\hfil {\footnotesize\draftcite{#1}}}\hskip 1cm
\if@filesw \immediate\write\@auxout
       {\string\bibcite{#1}{\the\value{\@listctr}}}\fi\ignorespaces}

\ifnum\extension=0 \newcommand\internal[1]{} \fi
\ifnum\extension=1 \newcommand\internal[1]{\mbox{}\\[-.25em]
                   {\small #1}\medskip} \fi

\def\citen#1{\if@filesw \immediate\write \@auxout {\string\citation{#1}}\fi%
\@tempcntb\m@ne \let\@h@ld\relax \def\@citea{}%
\@for \@citeb:=#1\do {\@ifundefined {b@\@citeb}%
    {\@h@ld\@citea\@tempcntb\m@ne{\bf ?}%
    \@warning {Citation `\@citeb ' on page \thepage \space undefined}}%
    {\@tempcnta\@tempcntb \advance\@tempcnta\@ne
    \setbox\z@\hbox\bgroup\ifcat0\csname b@\@citeb \endcsname \relax
    \egroup \@tempcntb\number\csname b@\@citeb \endcsname \relax
    \else \egroup \@tempcntb\m@ne \fi \ifnum\@tempcnta=\@tempcntb
    \ifx\@h@ld\relax \edef \@h@ld{\@citea\csname b@\@citeb\endcsname}%
    \else \edef\@h@ld{\hbox{--}\penalty\@highpenalty
    \csname b@\@citeb\endcsname}\fi
    \else \@h@ld\@citea\csname b@\@citeb \endcsname \let\@h@ld\relax \fi}%
\def\@citea{,\penalty\@highpenalty\hskip.13em plus.13em minus.13em}}\@h@ld}
\def\@citex[#1]#2{\@cite{\citen{#2}}{#1}}%
\def\@cite#1#2{\leavevmode\unskip\ifnum\lastpenalty=\z@\penalty\@highpenalty\fi%
  \ [{\multiply\@highpenalty 3 #1%
  \if@tempswa,\penalty\@highpenalty\ #2\fi}]}   %
\makeatother % end of CITE.STY

\newcommand\Ac[2]  {A^{#1,#2}_{}}

\def\aff           {affine Lie algebra}

\def\alphab        {{\bar\alpha}}
\def\alphaO        {\brev{\alpha}}
\newcommand\als[1] {\alpha^{(#1)}}

\def\AO            {\brev{A}}

\newcommand\AOc[2] {\brev A^{[#1],[#2]}_{}}

\def\atcha         {automorphism-twined character}

\def\avi           {a^{\Vee}_{i}}

\newcommand\avO[1] {\brev a^{\Vee}_{#1}}

\def\bearl         {\begin{array}{l}}
\def\bearll        {\begin{array}{ll}}
\def\bearlll       {\begin{array}{lll}}
\def\bearllll      {\begin{array}{llll}}

\def\BB            {{\hat B}}
\def\bll           {\mbox{$b_{\Lambda,\Lambda'}$}}
\def\bllch         {\mbox{$b^{(\ch)}_{\Lambda,\Lambda'}$}}

\def\bllo          {\mbox{$b^{[\omega]}_{\Lambda,\Lambda'}$}}

\def\bminus        {\mbox{${\cal B}_-$}}

\def\be            {\begin{equation}}

\def\bminus        {\mbox{${\cal B}_-$}}

\def\brev          {\breve}

\def\calh          {{\calH}}
\def\calH          {L}
\def\calw          {\mbox{$\Dot{\cal W}{.75}{.63}\;$}}
\def\calwc         {\mbox{${\cal W}_{\g,\g'}$}}
\def\calwoo        {\mbox{${\cal W}_{\g,\g'}^{(0)}$}}

\def\cc            {\Dot c{.49}{.31}}
\def\Ccs           {\Dot C{.51}{.41}}

\def\cd            {\Delta}

\def\cdo           {\brev\Delta}
\def\Cdot          {\!\cdot\!}
\def\cft           {conformal field theory}
\def\cfts          {conformal field theories}
\def\ch            {\Psi}
\def\Chi           {{\cal X}}
\def\CHI           {\Chi_{}^{[\omega]}}
\def\chii          {\raisebox{.15em}{$\chi$}}
\def\Chii          {{\cal X}}
\def\chil          {\chii_\Lambda^{[\omega]}}

\def\Chil          {\Chi_\Lambda^{[\omega]}}
\def\CHil          {\Chi_\lambda^{[\omega]}}

\def\Chili         {\Chi_{\lambda_i}^{[\omega]}}
\def\chilj         {\chii_{\lambda_j}^{[\omega]}}
\def\chillch       {\Chii_{\ilabl,\ch}^{}}

\def\chilP         {\chii_\LambdaP^{[\omP]}}
\def\ChilTau       {\overline\Chi^{[\omega]}(q)}
\def\ChilnTau      {\overline\Chi_{}^{[\omega]}(q^N)}
\def\chiO          {\raisebox{.15em}{$\brev\chi$}}
\def\ChiO          {\brev\Chi}

\def\chiv          {\chii_V^{}}
\def\cht           {{\tilde\Psi}}

\def\cO            {{\brev{c}}}
\def\CO            {{\bar C_2}}
\def\coo           {\mbox{$\Gamma_{\!00}$}}
\def\complex       {\mbox{$\dl C$}}
\def\compleX       {{\dl C}}

\def\CS            {{\cal S}}

\def\csa           {Cartan subalgebra}
\def\CT            {{\cal T}}

\def\cvira         {coset Virasoro algebra}

\def\deltaO        {{\brev\delta}}

\def\diag          {{\rm diag}}
\def\dO            {\brev{d}}
\def\DO            {\brev{D}}
\newcommand\Dot[3] {{#1}\hspace{-#3em}\raisebox{#2em}
                   {$\scriptscriptstyle\bullet$}}
\def\dpth          {{\rm dp}}

\def\dsum          {\displaystyle\sum}
\def\dyd           {Dynkin diagram}
\def\ee            {\end{equation}}
\def\eE            {{\rm e}}
\newcommand\EE[2] {E^{#1}_{#2}}
\def\eear          {\end{array}}
\def\Eins          {{\id}}
\def\eins          {{\scriptscriptstyle(1)}}
\def\el            {\ell}
\def\ela           {{\ell_{\bar\alpha}}}

\let\emb=\hookrightarrow
\def\emt           {energy-momentum tensor}
\def\Ene           {{\cal E}_1^{(\vec n)}}
\def\Enm           {{\cal E}_{}^{(\vec n,\vec m)}}
\def\Emz           {{\cal E}_2^{(\vec m)}}
\def\Ete           {\tilde{\cal E}_1^{(n_0,n_1,n')}}
\def\eps           {\epsilon}
\def\epsh          {\hat\epsilon}
\def\epsO          {\brev\epsilon}

\newcommand\erf[1]{(\ref{#1})}
\def\etaa          {{\eta_\alphab^{}}}
\def\facprop       {factorization property}
\newcommand\fc[3]  {{#1}\hspace{-#3em}\raisebox{#2em}{$\scriptstyle\circ$}}
\newcommand\fcft[3]{{{#1}^{\mskip-#3 mu\raise #2
pt\hbox{$\scriptstyle\circ$}}}}

\newcommand\fline  [1]{\vfill\noindent ------------------\\[1 mm]}
\def\fortria       {for $\tria\iN\{+,\circ,-\}$}
\newcommand\Frac[2]{\mbox{\large$\frac{#1}{#2}$}}
\def\fpr           {fixed point resolution}

\def\fssc          {\cite{fusS3}}
\def\fssn          {\cite{fusS4}}

\def\futnot#1      {\ifnum\draftcontrol=1
                   \footnote{~{\sc internal footnote:} #1}\ \fi}
\def\futnote#1     {\footnote{~#1}\ }
\def\g             {{\sf g}}

\def\gb            {\mbox{$\bar\g$}}
\def\gB            {{\bar\g}}
\def\gbP           {\mbox{$\gB'$}}

\newcommand\Gb[2]  {\bar G_{#1,#2}}
\def\gd            {\mbox{$\hat\g$}}
\def\gD            {\g_D^{}}

\def\geins         {\g_\circ^{(0)}}

\def\gH            {\g_\J^{}}
\def\Gid           {\mbox{$G_{\rm id}$}}
\def\GId           {{G_{\rm id}}}

\def\Gilabel       {G_{\Lambda,\Lambda'}}
\def\Gilabels      {G^\star_{\Lambda,\Lambda'}}
\def\Gjlabel       {G_{M,M'}}
\def\Gjlabels      {G^\star_{M,M'}}
\def\gK            {\g_K^{}}
\def\gkeins        {\g_K^{(0)}}

\def\go            {\g_\circ}
\def\goeins        {\g_\circ^{(0)}}
\def\gO            {\mbox{$\brev{\g}$}}
\def\GO            {\brev G}
\def\gOb           {\mbox{$\overline{\brev{\g}}$}}

\def\goh           {\hat\g_\circ}
\def\goheins       {\hat\g_\circ^{(0)}}

\def\gOo           {\brev\g_\circ}
\def\gOoh          {\hat{\brev\g}_\circ}

\def\gP            {\mbox{$\g'$}}

\def\gPO           {\brev\g'}
\def\gPP           {{\g'}}
\def\gpP           {\mbox{$\g'_+$}}
\def\Gstab         {{G_{\Lambda,\Lambda'}}}
\def\Gstabs        {{G^\star_{\Lambda,\Lambda'}}}
\newcommand\gstar[1] {\g_\circ^{\star (#1)}}
\def\gv            {\mbox{$g_{}^\Vee$}}
\def\gV            {g_{}^\Vee}

\def\gvO           {\mbox{$\brev g_{}^\Vee$}}
\def\gvo           {{\brev g^{\Vee}}}
\def\gz            {\mbox{$\g^{}_{\tria}$}}
\def\gzP           {\mbox{$\g'_{\tria}$}}
\def\h             {{\sf h}}
\def\half          {\mbox{$\frac12\,$}}
\def\hb            {\mbox{$\bar{\sf h}$}}

\newcommand\HH[2]{H^{#1}_{#2}}

\def\hill          {\mbox{${\calH}_\Lambda$}}

\def\hillO         {\mbox{${\calH}_\LambdaO$}}

\let\hl=\hill

\def\hll           {\mbox{${\calH}_{\Lambda,\Lambda'}$}}
\def\hllch         {\mbox{${\calH}^{(\ch)}_{\Lambda,\Lambda'}$}}
\def\hlo           {\mbox{${\calH}_{\om^\star\Lambda}$}}
\def\hlolo         {\mbox{${\calH}_{\omt\Lambda,\omtP\Lambda'}$}}
\def\hlP           {\mbox{${\calH}_{\el;\LambdaP}$}}
\def\hLP           {\mbox{${\calH}_{\LambdaP}$}}

\def\hlPt          {\mbox{${\calH}_{\taudd\el;\omdP\LambdaP}$}}

\def\hOO           {\brev{h}}
\def\HO            {\brev{H}}

\def\hsa           {horizontal subalgebra}
\newcommand\hsp[1] {\mbox{\hspace{#1 mm}}}
\def\hw            {highest weight}
\def\hwm           {highest weight module}

\def\hwv           {highest weight vector}
\def\hy            {$\mbox{-\hspace{-.66 mm}-}$}
\def\id            {\sl id}
\def\ii            {{\rm i}}
\def\ihwm          {irreducible highest weight module}

\def\I             {I}
\def\Ib            {\bar I}

\def\ilabel        {{(\Lambda;\Lambda')}}

\def\ilabl         {{[\Lambda;\Lambda']}}
\def\iN            {\!\in\!}

\def\IO            {{\brev{I}}}
\def\irmod         {irreducible module}

\def\J             {J}
\newcommand\JJ[1]  {\mbox{$\J^{#1}$}}

\def\jlabel        {{(M;M')}}

\def\jlabl         {{[M;M']}}

\def\kma           {Kac\hy Moo\-dy algebra}
\def\KO            {{\brev K}}

\def\kv            {\mbox{$k_{}^{\Vee}$}}
\def\kV            {k_{}^\Vee}
\def\KV            {k^\Vee}

\def\kvl           {\mbox{$k_\lambda^{\Vee}$}}

\def\kvO           {\mbox{$\brev k_{}^{\Vee}$}}
\def\kvo           {{\brev k^{\Vee}}}
\def\kvol          {{\brev k^{\Vee}_{\lambda}}}
\def\kvom          {{\brev k^{\Vee}_{\mu}}}
\newcommand\la[1]  {\mbox{$\Lambda_{(#1)}$}}
\newcommand\lab[1] {\mbox{$\Lambdab_{(#1)}$}}

\long\def\labl#1   {\label{#1}\ee \ifnum\draftcontrol=1
                   \mbox{ }\\[-12 mm]\query{#1}\\[5 mm] \fi}

\newcommand\lambdab{\bar\lambda}
\newcommand\Lambdab{\bar\Lambda}

\def\LambdaO       {{\brev\Lambda}}
\def\lambdaO       {{\brev\lambda}}

\def\lambdao       {\lambdaO}
\def\lambdaob      {\bar{\lambdaO}}
\def\Lambdaob      {\bar{\LambdaO}}

\def\lambdaOP      {{\brev\Lambda{}'}}
\def\lambdaP       {{\lambda'}}
\def\LambdaP       {{\Lambda'}}
\def\Lc            {\Dot L{.74}{.48}\hsp{.20}}
\def\Lcs           {\Dot L{.51}{.39}\hsp{.1}}

\def\li            {\mbox{$\Lambda_{(i)}$}}
\def\lie           {Lie algebra}
\def\Lie           {Lie group}
\def\llb           {\mbox{\large[}}
\def\Llb           {\mbox{\Large[}}

\def\LO            {{\brev{L}}}

\def\lP            {{\Lambda'}}
\def\LP            {L'}
\def\lrb           {\mbox{\large]}}
\def\Lrb           {\mbox{\Large]}}

\newcommand\mh[2]  {\hat m_{[#1],[#2]}}
\def\Mid           {\!\mid\!}
\def\mlambda       {m_\lambda^{[\omega]}}
\newcommand\mO[2]  {\brev m_{[#1],[#2]}}

\def\modinv        {modular invarian}
\def\Mu            {M}
\def\mub           {\bar \mu}

\def\muO           {\brev{\mu}}
\def\muo           {\brev{\mu}}
\def\MuO           {{\brev{M}}}
\def\muob          {\bar{\muo}}

\def\muOP          {{\brev{M}'}}
\def\muP           {\mu'}
\def\MuP           {M'}

\def\natnum        {\mbox{$\dl N$}}

\def\nextcase      {\\[.2em]\raisebox{.15em}{\rule{.4em}{.4em}}\, }
\def\nextcasE      {\vskip.1em\noindent\raisebox{.15em}{\rule{.4em}{.4em}}\, }
\def\nextcomm      {\vskip.2em\noindent\mbox{\raisebox{.32em}{\rule{.25em}
                   {.25em}}\hsp{-.34}\raisebox{.06em}{\rule{.25em}{.25em}}\, }}

\def\ocha          {twining character}
\def\ocond         {linking condition}

\def\olie          {orbit Lie algebra}
\def\om            {\omega}

\let\omchar=\ocha
\newcommand\omd[1] {{\dot\omega #1}}
\def\omD           {\dot\omega}
\newcommand\omdd[2]{{\dot\omega^{#1} #2}}

\def\omdi          {{\dot\omega i}}
\def\omdj          {{\dot\omega j}}

\newcommand\omdP[1]{{\omP{}^\star#1}}

\def\ome           {\omega_1}
\def\omh           {\omega_\h}

\def\omP           {{\omega'}}
\def\omT           {\mbox{$\omega^\star$}}
\newcommand\omt[1] {{\omega^\star #1}}
\def\oMT           {{\omega^\star}}
\newcommand\omtb   {\bar\omega^\star}

\def\omtLa         {{\omt\Lambda}}
\newcommand\omtP[1]{{{\omega'}^\star #1}}
\def\omz           {\omega_2}
\def\omzD          {\omega_2^\star}
\def\omzDP         {{\omega_2'}_{}^\star}
\def\one           {\mbox{\small $1\!\!$}1}
\def\onedim        {one-dimen\-sional}
\def\onehalf       {\mbox{$\frac12$}}

\def\onetonmr      {1,2,...\,,n-r}

\def\onetor        {1,2,...\,,r}
\def\otor          {0,1,...\,,r}

\def\otoNm         {0,1,...\,,N-1}
\let\p=\ell
\def\pbw           {Poin\-ca\-r\'e\hy Birk\-hoff\hy Witt}

\def\Pro           {\mbox{\sc p}^{}_{\!\omega}}

\def\proj          {\mbox{\sc p}^{\star-1}_{\!\omega}}
\def\projm         {\mbox{\sc p}^{\star}_{\!\omega}}

\def\prow          {\mbox{\sc p}_W}
\def\qc            {q_{}^{\Lcs_0-\Ccs/24}}

\long\def\query#1{\hskip 0pt{\vadjust{\everypar={}\small\vtop to 0pt{\hbox{}%
     \vskip -13pt\rlap{\hbox to 50.0pc{\hfil{\vtop{\hsize=8pc\tolerance=6000%
     \hfuzz=.5pc\rightskip=0pt plus 3em\noindent#1}}}}\vss}}}}%
\newcommand\rank[1] {\mbox{rank}\,#1}

\def\rep           {representation}
\def\Rep           {Representation}
\def\resp          {respectively}
\newcommand\restr[1] {|\raisebox{-.5em}{$#1$}}
\newcommand\restR[1] {|\raisebox{-.3em}{$#1$}}
\newcommand\rhob   {\bar \rho}
\def\rhoO          {\brev{\rho}}
\def\rhoob         {\,\bar{\!\rhoO}}
\def\rhs           {right hand side}

\newcommand\Scc[4] {{\cal S}^{}_{([#1;#2],\ch),([#3;#4],\cht)}}
\newcommand\sect[1] {\section{#1}\setcounter{equation}{0}}
\newcommand\Sect[2]{\sect{#1}\label{s.#2} \ifnum\draftcontrol=1
\query{s.#2}\fi}

\newcommand\Sf[4]  {\fc{\cal S}{.71}{.45}_{(#1;#2),(#3;#4)}}
\def\sgP           {\mbox{\small$\g'$}}
\newcommand\Skl[4] {{\cal S}^{}_{([#1;#2],\ch_k),([#3;#4],\ch_l)}}

\def\slz           {{SL(2,$\dl Z$)}}
\newcommand\Sm[2]  {S^{}_{#1,#2}}

\def\smat          {$S$-matrix}
\def\SNb           {{\cal S}^{[\omega^n]}}
\newcommand\Snb[4] {{\cal S}^{[\omega^n]}_{(#1;#2),(#3;#4)}}
\newcommand\Sno[4] {{\cal S}^{[\omega^0]}_{(#1;#2),(#3;#4)}}
\def\sO            {\brev s}

\newcommand\So[2]  {S^{[\omega]}_{#1,#2}}

\newcommand\Sob[4] {{\cal S}^{[\omega]}_{(#1;#2),(#3;#4)}}

\newcommand\Soo[4] {{\cal S}^{}_{(#1;#2),(#3;#4)}}
\newcommand\SoO[2] {\brev S^{}_{#1,#2}}
\newcommand\SoOP[2]{\brev S'_{#1,#2}}
\newcommand\SoP[2] {S^{[\omP]}_{#1,#2}}
\newcommand\SP[2]  {S'_{#1,#2}}

\def\sumgilabel    {\sum_{\om\in G_{\Lambda,\Lambda'}}}

\def\sumgilabels   {\sum_{\ch\in G^\star_{\Lambda,\Lambda'}}}

\newcommand\sumIO[1] {{\sum_{#1\in\IO}}}

\newcommand\sumrO[1] {\sum_{#1=0}^r}

\def\tauc          {\mbox{$\Dot\tau{.49}{.46}_\omega$}}

\newcommand\taudd[1]{\Dot\tau{.34}{.39}_\omega #1}

\def\tauo          {{\tau_\omega}}
\def\tauoP         {{\tau_\omP}}
\def\tbf           {twining branching function}
\let\tcha=\ocha
\def\Tcha          {Twining character}
\def\tchi          {\raisebox{.15em}{$\tilde\chi$}}
\def\tchil         {\tchi_\Lambda^{[\omega]}}
\def\ti            {{\brev I}}
\def\tildeR        {R^{\omega}_{}}
\def\tildER        {R^{\omega}}

\def\trhll         {{\rm tr}^{}_{{\calH}_{\Lambda,\Lambda'}}}

\def\tria          {\#}

\def\ttab          {{\bar\theta}}

\def\thetab        {\ttab}

\def\twodim        {two-di\-men\-si\-o\-nal}
\def\U             {{\sf U}}
\def\uaff          {untwisted affine Lie algebra}
\def\ug            {\mbox{${\sf U}(\g)$}}
\def\ugtw          {\mbox{$\tilde{\sf U}(\g)$}}
\def\untw          {^{(1)}}

\def\Vee           {{\scriptscriptstyle\vee}}
\newcommand\version[1] {\ifnum\draftcontrol=1 \typeout{}\typeout{#1}\typeout{}
                   \vskip3mm \centerline{\fbox{\tt DRAFT -- #1 -- \today}}
                   \vskip3mm \fi}

\def\vlaP          {\mbox{$v^{}_{\el;\lambdaP}$}}

\def\vmL           {{V_\Lambda}}

\def\vml           {\mbox{$V_\Lambda$}}
\def\vmlt          {\mbox{$V^\omega_\Lambda$}}

\def\wh            {\hat w}
\def\Wh            {\hat W}

\def\whi           {\hat w_{[i]}}
\def\whj           {\hat w_{[j]}}
\def\Wl            {W_{(\lambda)}}

\def\wi            {w_i}
\def\womi          {w_{\omD i}}
\def\wO            {\brev{w}}
\def\wOi           {\brev{w}_{[i]}}
\def\wOj           {\brev{w}_{[j]}}
\def\WO            {\brev{W}}
\def\wrt           {with respect to }
\def\wrtt          {with respect to the }
\def\WZW           {Wess\hy Zumino\hy Witten}

\def\wzwt          {WZW theory}
\def\wzwts         {WZW theories}

\def\zet           {{\dl Z}}

\def\zetplus       {\mbox{${\zet}_{>0}$}}
\def\zetpluso      {\mbox{${\zet}_{\geq 0}$}}
\def\zwei          {{\scriptscriptstyle(2)}}

\global\def\draftcontrol{0}
\catcode`\@=12

\documentstyle[12pt]{article}

\begin{document}
%%% for draft versions, suppress in definitive version
%\draft
\version\versionno

\let\dl=\bf
\input amssym.def \input amssym.tex \def\dl{\large\Bbb }
% if you do not have these files, or the fonts used there, remove line above

%%%%%%%%%%%%%%%%%%%%%%%%%%%%%%%%%%%%%%%%%%%%%%%%%%%%%%%%%%%%%%%%%%%%%%%
\begin{flushright}  {~} \\[-23 mm] {\sf q-alg/9511026} \\
{\sf DESY 95-222}\\{\sf NIKHEF/95-062}\\{\sf IHES/P/95/97}
\\[1 mm]{\sf November 1995} \end{flushright} \vskip 2mm

\begin{center} \vskip 11mm
{\Large\bf
TWINING CHARACTERS, ORBIT LIE ALGEBRAS,} \vskip3mm
{\Large\bf \ AND FIXED POINT RESOLUTION $^{\sf T}$}
\vskip 15mm
{\large J\"urgen Fuchs, $^{1,\sf H}$ \ \,
 Bert Schellekens, $^2$\, \ Christoph Schweigert $^3$} \\[5mm]
{$^1$ \footnotesize DESY, Notkestra\ss e 85, \ D -- 22603~~Hamburg} \\[.4em]
{\hsp{2.5}$^2$ \footnotesize NIKHEF-H, Kruislaan 409,\hsp{.3} \
NL -- 1098 SJ~~Amsterdam }\\[.4em]
{\hsp{3.5}$^3$ \footnotesize IHES, 35, route de Chartres, \
F -- 92440~~Bures-sur-Yvette} \end{center} \vskip 22mm
\begin{quote}{\bf Abstract}.\\
We describe the resolution of field identification fixed points in coset
conformal field theories in terms of representation spaces of the coset
chiral algebra. A necessary ingredient from the representation theory of
Kac\hy Moody algebras is the recently developed theory of twining characters
and orbit Lie algebras, as applied to automorphisms representing identification
currents.
\end{quote}
%\vskip 9mm
\vfill {}\fline{} {\small
$^{\sf T}$~~Based on lectures by J.\ Fuchs at the workshop `New Trends in
Quantum Field Theory'\\{}\hsp{3.4} (Razlog, Bulgaria, August 1995)
\\[.4em] $^{\sf H}$~~Heisenberg fellow} \newpage

\section{Fixed point resolution}

While in the title of this paper we refer to some specific mathematical
structures, namely certain (recently developed)
aspects of the \rep\ theory of Kac\hy Moody \lie s,
the problem that we address is a rather general one, and in fact the basic
ideas can be described without going into any technical details.
This problem, to which we refer as {\em \fpr}, arises whenever
one has to mod out redundancies in the description of a system
(`gauge symmetries') in a situation
where the orbits of the redundancy transformations have unequal
sizes, and hence potentially arises in various different areas of
physics and mathematics.

The presence of a redundancy symmetry
means that distinct combinations of the, a priori basic, data of a
theory actually describe the same physical situation. Thus the physical
state of the system in question is not described by an individual set of
those data, but rather by an equivalence class of such sets that is given by
an orbit of the redundancy symmetries. For instance, in Yang\hy Mills
theory physical states do not correspond to individual configurations
of the fundamental fields, but rather to gauge orbits thereof.
In short, the prescription is
  \be  \mbox{\fbox{\ individual configurations $\Phi$\ }} \quad \longrightarrow
  \quad \mbox{\fbox{\ orbits\hsp{.4} \raisebox{.15em}{$\scriptstyle[$}\,%
  $\Phi$\,\raisebox{.15em}{$\scriptstyle]$}\ }} \hsp3  \labl{fig1}
The prescription \erf{fig1} must, however, possibly be modified
as soon as orbits of different sizes are present.
Such a situation occurs e.g.\ in Yang\hy Mills theory, where
reducible connections lead to `shorter' orbits, implying that the
space of gauge orbits is not a smooth manifold, but a stratified variety
which has singularities. While a priori it may well be consistent to work
with the space of orbits even in this more complicated situation,
one should not be surprised if doing so one encounters problems;
e.g.\ in the Yang\hy Mills case one may think of
ambiguities in the quantization procedure.

Here we consider a situation where the dynamics is under control so that
the inconsistency of a naive implementation of the redundancies can be
easily detected.
Namely, we analyze coset theories, i.e.\ \cfts\ whose Virasoro algebra is
obtained by subtracting from the Sugawara \emt\ of a \wzwt\ based on an affine
\lie\ \g\ the Sugawara tensor of the \wzwt\ based on a subalgebra \gP\ of \g.
The redundancy symmetry in this case includes, but is not exhausted by,
the subalgebra \gP.
In the presence of orbits of different sizes, in these theories one would
not obtain a modular invariant partition function
when proceeding according to the prescription \erf{fig1}. One advantage of the
type of theory we consider is that the presence of different orbit sizes can
be attributed to the action of a discrete finite group.

The resolution of this problem is to supplement the recipe \erf{fig1} by a
`resolution of the fixed points'.
By this we simply mean that the states of the theory must be organized
in a more complicated way than suggested by \erf{fig1}.
This somewhat vague statement can be made explicit
in the \cft\ context, where it means that some of the naively obtained
modules of the coset chiral algebra are {\em not\/} irreducible.
Rather, the irreducible modules are submodules which are eigenspaces of the
additional redundancy transformations. For the full theory obtained by
combining its
two chiral halves this implies that not all states one would naively expect
are present, i.e.\ that on the space of states an additional projection
takes place.

We defer a more precise description of these issues to section \ref{s.cos}.
Before that, we have to describe certain mathematical structures which are
basic ingredients for the analysis of section \ref{s.cos}. The
key concepts are {\em \tcha s\/} and {\em \olie s\/} which have been
introduced in \fssc; these structures are also interesting in their own right.

\sect{\kma s and diagram automorphisms}

Both \tcha s and \olie s arise in a much broader context than the one of
affine \lie s that is relevant to the application in \cft;
they are both determined by two basic data:\\
(1) a symmetrizable \kma\ \g, \ and \
(2) a diagram automorphism $\om$ of \g.  \smallskip

A {\em symmetrizable \kma\/} is a \lie\ which possesses
both a Cartan matrix and a Killing form. That is, there is a square matrix
$A$, the Cartan matrix of \g, which has diagonal entries $\Ac ii=2$ and
non-positive integral off-diagonal entries such that $\Ac ij=0$ iff $\Ac ji=0$,
as well as an invertible diagonal matrix $D$ such that $DA$ is symmetric.
The algebra \g\ associated to $A$ is obtained as follows \cite{KAc3}.
The \csa\ $\go$ is by definition an abelian Lie algebra of
dimension $2n-r$, where $n$ is the dimension and $r$ the rank of $A$;
there are (uniquely up to isomorphism)
$n$ linearly independent elements $H^i$ of $\go$ and
$n$ linearly independent functionals $\als i\iN\go^\star$ (the simple
roots) such that $\als i(H^j)=\Ac ij$ for $i,j=1,2,...\,,n$;
\g\ is generated freely by $\go$ and by elements $\EE i\pm\equiv
E^{\pm\alpha^{(i)}}$, with $i\iN\{1,2,...\,,n\}$, modulo the relations
  \be \begin{array}{ll}
  [x,y] =0  \quad\mbox{for all}\ x,y\in\go \,,\ \ &
  [x,\EE i\pm]=\pm\als i(x)\,\EE i\pm \quad\mbox{for all}\ x\in\go \,,
\\[2.5mm]
  [\EE i+,\EE j-]=\delta_{ij}\, \HH j{} \,, &
  {({\rm ad}^{}_{\EE i\pm})}^{1-\Ac ji}_{}\EE j\pm =0 \quad{\rm for}\ i\ne j
  \,. \end{array}\labl C
The \kma\ \g\ has a triangular decomposition $\g = \g_+ \oplus \go \oplus
\g_-$,
where $\g_\pm$ are subalgebras generated freely by the $\EE i\pm$
modulo the Serre relations, i.e.\ the last type of relations in \erf C.
For any null eigenvector of the Cartan matrix, the corresponding linear
combination of the $H^i$ is a central element, and hence the central subalgebra
$\gK$ of \g\ has dimension $n-r$.
The derived algebra $\gd:=[\g,\g]$ of \g\ contains $\gK$, and it has
a triangular decomposition, namely $\gd= \g_+ \oplus \goh \oplus \g_-$,
where $\goh$ is the span of the elements $H^i$, $i=1,2,...\,,n$.
The generators of a complement $\gD$ of $\goh$ in $\go$ are called the
derivations of \g.

The \dyd\ of \g\ is the graph with $n$ vertices that has coincidence matrix
$2\cdot\one-DA$.  Without loss of generality, we will assume that
the \dyd\ is connected, i.e.\ that $A$ is indecomposable. A {\em diagram
automorphism\/} $\om$ of \g\ is an automorphism that
acts on the generators corresponding to simple roots like
  \be   \omega(E^i_\pm) :=  E^{\omd i}_\pm     \labl{oei}
by a symmetry of the \dyd, i.e.\ by a permutation $\omD$ satisfying
$\Ac\omdi\omdj=\Ac ij$ for all $i,j$.
The mapping \erf{oei} extends to automorphisms of $\g_+$ and $\g_-$, and by
$\omega(H^i):= \omega([E^i_+,E^i_-]) = [ E^{\omd i}_+, E^{\omd i}_-]
= H^{\omd i}$ it extends to a unique automorphism of \gd\ of finite order.

The extension of $\omega$ to the rest of \g\ is in general not unique.
To see this, we choose a basis of $n-r$ eigenvectors $K^a$, $a=\onetonmr$,
of the central subspace $\gK$ such that
  \be  \omega(K^a) = \zeta^{n_a} K^a  \qquad{\rm with}\quad
  \zeta:=\exp(2\pi\ii/N)  \labl{zeta}
($N$ is the order of $\omD$); this is possible because $\om$ maps
$\gK$ bijectively to itself. This extends to a basis
of $\goh$ with further generators $\JJ p$, $p=\onetor$, satisfying
$\omega(\JJ p) = \zeta^{m_p} \JJ p$. The restriction of the non-degenerate
invariant bilinear form $(\cdot\mid\cdot)$ on $\go$ to the span $\gH$ of the
\JJ p is non-degenerate, and $(K^a\mid x) =0$ for all $x\in\goh$.
It follows that there are $n-r$ unique elements $D^a$ of $\go$ such that
$(D^a\Mid K^a) = 1$ and $(D^a\Mid\cdot\,) = 0$ for all other basis elements
of $\go$. The $D^a$ form a basis for a complement $\gD$ of
$\goh$ in $\go$. Using the automorphism property of $\om$, we find
  \be  \omega(D^a)  = \zeta^{-n_a} D^a + \sum_{b=1}^{n-r} M^a_b \zeta^{n_b} K^b
  \,,\labl{ans2}
with $M$ an antisymmetric $(n-r)\times(n-r)$ matrix.

Thus the freedom in defining $\om$ consists in adding central terms to
$\omega(D^a)$ and is parametrized by the antisymmetric matrix $M$. In
particular, if \g\ is simple, affine, or hyperbolic, $\om$ is in fact
uniquely determined by $\omD$. Also, when acting on $\gd$,
$\om$ has the same order $N$ as $\omD$, and hence in the general case
we can (and will) restrict the freedom in $\om(D)$ further by requiring
that $\om$ has order $N$ on all of \g;
then $M^a_b$ must vanish whenever $n_a=-n_b\;{\rm mod}\,N$.
It is quite important that the only freedom in the action of $\om$
consists in adding central terms.
In the highest weight representations of our interest, the central elements act
as multiples of the identity; in particular, the character-like quantities we
will consider are affected by the freedom in choosing $\omega$
only via multiplication with an overall factor.

\sect{\Tcha s}

The ordinary character $\chiv$ of a \g-module $V$ on which the
action of the \csa\ can be diagonalized is a (formal) function on the
\csa\ $\go$, defined by
  \be  \chiv: \quad \go   \to  \complex \,,\quad
  \chiv(h):=   {\rm tr}_V^{} \eE^{2\pi\ii R_\Lambda(h)} \,.  \ee
($R$ is the \rep\ in which \g\ acts on $V$; from now on we suppress
the symbol $R$ whenever no confusion can arise.) The character can be expanded
in terms of the weights $\lambda$ of $V$ as $\chiv =\sum_\lambda m_\lambda^{}
\,\eE^{2\pi\ii\lambda}$, where $m_\lambda$ is the
multiplicity of $\lambda$, i.e.\ the dimension of the subspace $\Wl$
of $V$ that has weight $\lambda$.
Hence $\chiv$ is the generating functional for the dimensions of the weight
spaces. For \hwm s \vml\ with \hw\ $\Lambda$, $m_\lambda\ne0$ implies that
$\lambda\le\Lambda$, i.e.\ that $\Lambda-\lambda$ is a linear combination
of simple roots with non-negative integral coefficients.

Now the idea of the \tcha\ is to count not just multiplicities of states,
but in addition to keep track of the action of $\om$, analogous as for a
character-valued index. To put this idea to work, we must specify what
`the action of $\om$' on $V$ means. This map, to be
denoted by $\tauo$, is the natural induced action
on the \g-module $(V,R)$, which means that $\tauo$ does not change $V$ as an
abstract vector space while it does change the \rep\ $R\!:\ \g \to {\sl
End}(V)$
of \g\ by endomorphisms $R(x)$ of $V$; namely,
the \rep\ of \g\ on $V$ changes according to the map
  \be  R(x)\;\mapsto\; \tildeR(x) := R(\om(x))  \labl T
for all $x\in\g$. Thus in general the map does change the
(isomorphism class of the) module.
 \futnote{In other words, $R(\g)$ and $\tildeR(\g)$ describe two, generically
inequivalent, embeddings of \g\ into the algebra {\sl End}($V$). It is only
after applying to $(V,R)$ and $(V,\tildeR)$ the forgetful functor from the
category of \g-modules to the category of vector spaces
that $(V,R)$ and $(V,\tildeR)$ become identical objects.}
Further, the action of $\omega$
preserves the triangular decomposition of \g\ so that
$(V,\tildER_\Lambda)$ is again a Verma module if $(V,R_\Lambda)$ is.
Moreover, the sets of primitive singular vectors of $(V,R_\Lambda)$ and
$(V,\tildER_\Lambda)$ coincide, and hence acting on the
irreducible quotient of a Verma module yields again an irreducible
module. To save space, for the time being we will denote both the Verma module
with \hw\ $\Lambda$ and its irreducible quotient by the same symbol \vml.

While the \hwv\ in the modules $(V,R)$ and $(V,\tildeR)$ is one and
the same element $v^{}_{\rm h.w.}$ of the underlying vector space $V$,
its associated weight gets transformed by the dual action $\omT$ of $\om$
on the weight space $\go^\star$ of \g\ (this mapping $\beta\mapsto\omT\beta$
reads $(\omT\beta)(x)=\beta(\omega^{-1}x)$ for all $\beta\iN\go^\star$ and
all $x\iN\go$). Thus, while $v^{}_{\rm h.w.}$ has weight $\Lambda$ in
$(V,R_\Lambda)$, it has weight $\omT\Lambda$ in $(V,\tildER_\Lambda)$.
In other words, $(V, R_\Lambda)$ is isomorphic to the module \vml,
while $(V,\tildER_\Lambda)$ is isomorphic to $V_\omtLa$.
Via these isomorphisms one and the same vector $v\iN V$
is identified with an element $v'$ of \vml\ and another element $v''$ of
$V_\omtLa$. We therefore define $\tauo$ as the map
$\tauo\!:\,\vml\to\vmlt$ that acts as $v'\mapsto v''$; then $\tauo(R_\Lambda
(x)\cdot v)= R_\omtLa(\omega(x))\cdot\tauo(v)$, i.e.\ the diagram
  \be
  \begin{array}{rcl} \vml &
  \stackrel{R_\Lambda(x)}{\mbox{---------}\!\!\!\longrightarrow}
  & \vml \\[.5mm]
  {\scriptstyle \tauo}\,\downarrow\ && \ \downarrow\,{\scriptstyle \tauo}
  \\[.5mm]  V_\omtLa \! & \stackrel{R_\omtLa(\om(x))}
  {\mbox{---------}\!\!\!\longrightarrow} & \! V_\omtLa
  \end{array} \labl{C'}
commutes for all $x\in\g$ and all $v\in\vml$. Thus $\tauo$ behaves as a
generalization of an intertwiner map; accordingly we refer to \erf{C'}
as the {\em $\omega$-twining property\/} of $\tauo$.

The proper implementation of the idea to keep track of the action of $\om$
in the character is thus expressed by the following definition.
The {\em\tcha\/} (or \atcha) $\Chil$
of a (Verma or irreducible) module \vml\ is the (formal) function
  \be  \Chil: \quad \go \to \complex \,,\qquad
  \Chil(h):=  {\rm tr}_\vmL^{} \tauo\, \eE^{2\pi\ii R_\Lambda(h)} \labl{Chil}
on $\go$. Analogously as for the ordinary character $\Chi$, we can write
the \omchar\ as
  \be \Chil = \sum_{\lambda\leq\Lambda} \mlambda\, \eE^{2\pi\ii\lambda} \,.
  \labl{weide}
Here $\mlambda$ is non-zero only if $\omega^\star(\lambda)=\lambda$, in
which case it is the trace of the restriction of $\tauo$ to the
weight space $\Wl$ of weight $\lambda$;
thus the twining character is the generating functional for these traces.
Strictly speaking, the notation `tr' makes sense only if
$\omT\Lambda=\Lambda$, since only in this case
$\tauo$ is an endomorphism of a single \hwm\ rather than
a map between two different ones. We will refer to such \g-weights $\Lambda$
as {\em symmetric\/} weights, and simply define $\Chil$ to be zero for
non-symmetric weights; below we will be mainly interested in
symmetric weights. Generically some contributions to $\Chil$ have non-zero
phase, so that is not at all obvious that
the term character is an appropriate name for the object defined by
\erf{Chil}. However, as we will see below, under suitable conditions
the expansion coefficients $\mlambda$ are non-negative integers, so that
the term character is indeed justified.

The following properties of \tcha s can be derived in a straightforward manner.
\nextcase
The \tcha\ is majorized by the ordinary character of the same module, and hence
in particular its domain of convergence contains the one of the ordinary
character.
\nextcase
Since $\mlambda\ne0$ implies $\omT\lambda=\lambda$,
we can restrict the sum in \erf{weide} to symmetric weights.
\nextcase
Together with the cyclic invariance of the trace, the $\om$-twining
property \erf{C'} implies that
  \be  \Chil(h) = \Chil(\om(h)) \labl{39}
for all $h\iN\go$.
\nextcase
Using the eigenspace decompositions
 $\go = \bigoplus \go^{(j)}$ and $\go^\star = \bigoplus \gstar j$, with
 $\om x =\zeta^j x  \;\forall\;x\iN\go^{(j)}, \;
  \omT\beta=\zeta^j\beta \;\forall\;\beta\iN \gstar j$,
it follows that the symmetric weights, i.e.\ those in $\gstar0$, are non-zero
only on the symmetric part $\go^{(0)}$ of the \csa. This implies that
$\Chil(h) = \sum_{\lambda} \mlambda \,\eE^{2\pi\ii\lambda}(\sum_j h^{(j)})
  = \sum_\lambda \mlambda\,\exp[ 2\pi\ii\,\lambda(h^{(0)})]$,
so that \erf{39} can be strengthened to
  \be  \Chil(h)=\Chil(h^{(0)}) \,.  \labl{chilh0}
Thus the \omchar\ depends on $h\iN\go$ non-trivially only through
its component $h^{(0)}$ in $\go^{(0)}$.
Correspondingly, from now on we will view $\Chil$ no longer as
a function on $\go$, but on $\go^{(0)}$.

\sect{Orbit \lie s}

To analyze the \omchar s in more detail, we need another concept, namely
the notion of an {\em \olie}. This concept is based on the simple observation
that any symmetry $\omD$ of a \dyd\ divides the set of nodes of the diagram
into invariant subsets of size $N_i$, the orbits $[i]$ of $\omD$.
Now these orbits can be viewed as the nodes of another \dyd, obtained by an
appropriate `folding'. The correct prescription includes suitable weight
factors; in terms of the associated Cartan matrix $A$, it reads as follows:
\nextcase
sum up the $N_i$ rows of $A$ belonging to each orbit $[i]$;
\nextcase
multiply the resulting row
with the integer $\;s_i:=1-\sum_{l=1}^{N_i-1}\!\Ac{\omdd li}i\in\zetplus$\,;
\nextcase
eliminate redundant columns.\\[.3em]
Thus we obtain a matrix $\AO$ whose rows and columns are labelled by the
$\omD$-orbits,
with entries
  \be  \AOc ij := s_i \frac{N_i}N \sum_{l=0}^{N-1}\Ac{\omdd li}j \,. \labl{AO}
In general, this is no longer a symmetrizable Cartan matrix. However, we can
prove that $\AO$ is again a symmetrizable Cartan matrix whenever
  \be  s_i \in \{1,2\} \quad {\rm for\ all}\ i \,;\labl{ocond}
this requirement will be called the {\em \ocond}.
In terms of the \dyd, \erf{ocond} means that each of
the nodes on an orbit is either connected by a single link to precisely
one node on the same orbit or not linked at all to other nodes on the same
orbit. In the sequel we will often restrict our attention to the class of those
\dyd\ symmetries $\omD$ which do satisfy the \ocond.

\internal{
The proof that condition \erf{ocond} implies that $\AO$ is a symmetrizable
Cartan  matrix is not difficult. For instance, one has immediately
$\AOc ii = s_i\,(3-s_i) = 2$, and the integrality of the entries of $\AO$
is made manifest by rewriting \erf{AO} as $\AOc ij = s_i\sum_{l=0}^{N_i-1}
\Ac{\omdd li}j$. Also, if $D=\diag(d_i)$ is
a non-singular diagonal matrix such that $DA$ is symmetric, then
$\DO = \diag(\dO_{[i]})$ with $\dO_{[i]}:=Nd_i/s_iN_i$ is a
non-singular diagonal matrix such that $\DO\AO$ is symmetric.
}

Symmetrizable Cartan matrices are either of {\em finite\/} (or {\em simple})
type (then the matrix $DA$ is positive definite), of affine type
(then $DA$ is positive semidefinite and has exactly one eigenvector with
eigenvalue zero), or else of {\em indefinite\/} type.
One can see that the matrix $\AO$ obtained from $A$ by the prescription
\erf{AO} is always of the same type as $A$. Also, by inspection one finds that
all symmetries of simple and affine \dyd s satisfy the \ocond,
except for the order $N$ automorphisms of the \aff s $A_{N-1}\untw$ which
rotate the \dyd\ (but the latter will also be dealt with
separately below). Further, for the indefinite Cartan matrices which
are {\em hyperbolic\/} (i.e., characterized by the property that any
connected subdiagram of the \dyd\ that one obtains by deleting any node from
the \dyd\ of \g\ is finite or affine), again also $\AO$ is hyperbolic.

We denote the \kma\ corresponding to the Cartan matrix $\AO$ by \gO\ and call
it the {\em orbit \lie} associated to \g\ and $\om$. A crucial observation is
that the \csa\ $\gOo$ of \gO\ can be related to the subalgebra $\go^{(0)}$
of the \csa\ $\go$ of \g\ that is invariant under $\omega$ by defining, for
$h=\sum_i v_i H^i\iN\goheins\equiv\geins\cap\goh$,
  \be  \Pro(h) := \sum_{[i]} N_i v_i\, \HO^{[i]}  \, \labl{Pro}
(this does not depend on the choice of representatives of the orbits $[i]$,
because $v_i = v_{\omD^l i}$ for $h\iN\goheins$).
The map $\Pro$ has the following properties:
\nextcase
The mapping\, $\Pro\!:\; \goheins\to\gOoh$\, is one-to-one.
\nextcase
The invariant bilinear form on $\gOoh$ is proportional to the restriction
to $\goheins$ of the invariant bilinear form on $\goh$:
for all $h_1,h_2\iN \goheins$ we have
  \be   (\Pro(h_1)\mid \Pro(h_2)) = N\,(h_1\mid h_2) \,. \labl{normrel}
\nextcasE
If $K\iN\goheins$ is central, then so is $\Pro(K)\iN\gOoh$, and vice versa.
Thus in particular the dimension of the space $\gkeins$ of invariant
central elements is the number of $\omD$-orbits minus the rank of $\AO$.
\nextcase
As a consequence, one can employ the non-degeneracy of the invariant
bilinear form on $\goeins$ to extend the map $\Pro$ to all of $\go$
in such a way that \erf{normrel} is still valid.
\nextcase
By duality, the correspondence between $\go^{(0)}$ and $\gOo$
implies analogous relations for the weight spaces. Thus there is
a bijective map\, $\projm\!\!:\,\gOo^\star\to \gstar0$\,
between the weights of $\gO$ and the symmetric weights $\lambda\iN\gstar0$,
such that, analogously to \erf{normrel},\, $(\lambda\mid\mu) = N\Cdot
(\proj(\lambda)\mid \proj(\mu))$.
For brevity, we will often denote the pre-image\, $\proj(\lambda)\iN
\gOo^\star$\, of a symmetric weight $\lambda$ simply by $\lambdaO$.

\sect{\Tcha s as characters of \olie s}

The two structures introduced above -- \tcha s and \olie s --
are both determined through a symmetrizable \kma\ \g\ and a diagram
automorphism $\om$ of \g, so that one should not be surprised that
they are closely related. Nevertheless we think that
the precise form of the relationship is
most remarkable. Namely, provided that $\om$ satisfies the \ocond, we find,
both for Verma modules and for their irreducible quotient modules: \vskip.5em

{}\hsp5
\fbox{\parbox{34em}{The \tcha s of the \hwm s over \g\ coincide with the
ordinary characters of the \hwm s over the \olie\ \gO.}}

\vskip.7em \noindent
To specify what `coincide' means, the relevant `projections' needed for this
statement to make sense must be performed, i.e.\ more explicitly the
assertion of this theorem is that
  \be
  \mbox{\fbox{$\quad
  \Chil(h) = \ChiO_{\LambdaO}(\hOO)   \qquad {\rm with} \qquad
  \LambdaO=\proj(\Lambda) \quad{\rm and}\quad \hOO=\Pro(h) \,.\quad $}}
  \labl{1a}
An immediate consequence of this identification is that the term
`character' for $\Chil$ is indeed justified; also, under suitable
circumstances $\Chil$ has nice modular transformation properties.

As already mentioned, the \ocond\ is not satisfied for the
order $N$ diagram automorphisms of $A_{N-1}^{(1)}$ which do appear in the
application to \cft.  But we can treat these cases as well, and the result is
again very simple: the only non-zero term in the expansion \erf{weide} for
both the irreducible and Verma modules is the one for the \hw, i.e.\ we have
  \be  \mlambda=0 \quad{\rm for}\ \lambda\ne\Lambda \,. \labl2

The proof of the statements \erf{1a} and \erf2 is a mixture of a
`conceptual proof' and of `verification' by explicit calculations. We will
concentrate on the conceptual parts; the computational details can be found
in \cite{fusS3}. A crucial ingredient is the identification of a natural
action of $\WO$, the Weyl group of the orbit \lie\ $\gO$, on the \omchar s.
\nextcase
We claim that the Weyl group of $\gO$ is isomorphic to a subgroup of the Weyl
group of \g, and that the action of this subgroup on $\go^\star$ commutes
with the action of $\omega^\star$.  \\[.2em]
We denote by $\wi$ and $\wOi$ the fundamental
reflections which generate the Weyl groups $W$ of \g\ and $\WO$ of \gO, \resp,
and employ the map $\projm$ to push the action of $\WO$ on
the weight space $\gOo^\star$ of the \olie\ \gO\ to an action on $\gstar{0}$.
First we show that for any $\wOi\iN\WO$ there exists an associated
element $\whi$ of $\Wh$
which acts on $\gstar 0$ precisely like $\wOi$ acts on $\gOo^\star$, i.e.\
  $\proj(\whi(\lambda)) = \wOi(\proj(\lambda))$
for all $\lambda\iN\gstar0$. Now by direct calculation, we obtain
$\projm(\alphaO^{(i)^\Vee}) = \frac N{N_i} \sum_{l=0}^{N_i-1}
\alpha^{(\omD^l i)^\Vee}$ and
$\projm(\alphaO^{(i)}) = s_i\sum_{l=0}^{N_i-1} \alpha^{(\omD^l i)}$, which
implies that the latter relation is equivalent to requiring that the action
of $\whi$ reads
  \be  \whi(\lambda) = \lambda - s_i\cdot\sum_{l=0}^{N_i-1}
  (\lambda\mid \alpha^{(\omD^l i)^\Vee})\, \alpha^{(\omD^l i)}  \, . \labl{w12}
We denote the mapping $\wOi\mapsto\whi$ by $\prow$.
As can be checked by employing the relations in the Weyl group $W$, the ansatz
  \be  \whi = \left\{ \bearll
  \prod_{l=0}^{N_i-1} w_{\omD^l i} & {\rm for}\ s_i=1 \,, \\[.2em]
  \wi\,\womi\,\wi = \womi\,\wi\,\womi & {\rm for}\ N_i=s_i=2 \,
  \eear\right.\labl{whi}
does fulfill \erf{w12}.\,\,%
 \futnote{
For $s_i=2$, there is precisely one $m\iN\{1,...\,, N_i-1\}$ with
$\Ac{\omdd mi}i = -1$. It follows that $\omdd mi=\omdd {-m}i$, which in
turn implies that the orbit length $N_i$ is even and that
the restriction of the \dyd\ of \g\ to this orbit is the \dyd\ of
$N_i/2$ copies of $A_2$. As a consequence, without loss of
generality we can restrict ourselves to the case $N_i=2$. Otherwise we first
treat the automorphism $\omD^{N_i/2}$, which has order two and possesses
$N_i/2$ orbits each of which corresponds to the \dyd\ of $A_2$. On the
set of orbits of $\omD^{N_i/2}$,
the automorphism $\omD$ induces an automorphism $\ddot\omega$
of order $N_i/2$; all orbits \wrt $\ddot\omega$ have $s_j=1$.}
(The relations of $W$ also imply that for $s_i=1$
the multiple product in \erf{whi} does not depend on the order of the factors.)

We can also check that the $\whi$ are reflections, i.e.\ square to the
identity.
Moreover, $\whi$ commutes with $\omega^\star$, $[\whi,\omega^\star]=0$
for all $i$; this implies in particular that the action of $\whi$ respects the
orbits of $\omega^\star$.
Next we define $\Wh$ as the subgroup of $W$ that is generated
by the elements $\whi$. To show that the map $\wOi\mapsto\whi$
extends to a  group isomorphism $\prow\!:\;\WO\to\Wh$,
we view the Weyl group $\WO$ as the Coxeter group that is freely generated
by the generators $\wOi$ modulo the relations $(\wOi)^2 = \Eins$ and
  $(\wOi\wOj)^{\mO ij} = \Eins$ for all $i,j$ with $i\neq j$,
where $\mO ij= 2,3,4,6$ for $\AOc ij\AOc ji = 0,1,2,3$ and (formally)
$\mO ij=\infty$ for $\AOc ij\AOc ji\ge4$. We can prove that $\Wh$ obeys
identical relations, i.e. is a Coxeter group with $\mh ij=\mO ij$,
and hence $\Wh\cong\WO$.

\internal{
What is in fact easy to see is that $\mh ij\ge\mO ij$. Namely,
if $(\whi\whj)^{\mh ij}=\Eins\iN W$, then in particular $(\whi\whj)^{\mh ij}$
acts as the identity on $\gstar0$; then
also $(\wOi\wOj)^{\mh ij}\iN\WO$ acts as the identity on the
weight space of $\gO$, hence is the identity element of $\WO$;
this implies that $\mO ij$ must be a divisor of $\mh ij$.
The inequality $\mh ij\ge\mO ij$ already proves the assertion
for $\AOc ij\AOc ji\geq 4$. In the remaining cases one can show,
by exploiting the explicit form of the Coxeter
relations for $W$, that also $\mh ij\le\mO ij$, so that equality
follows again. However, this necessitates a lengthy case by case study
depending
on the value of $\Ac ij\Ac ji$ (and also of $s_i$, $s_j$, $N_i$, $N_j$)
\cite{fusS3}; the details can be found in \cite{fusS3}.
}
\nextcase
The next step consists in considering the action of the group $\Wh$ on
the \tcha s\ that is induced via the action \erf{w12} on \g-weights. \\[.2em]
{}From this point on the proof proceeds in a way rather similar to
Kac' proof of the Weyl\hy Kac character
formula, see e.g.\ \cite[pp.\,152,\,172]{KAc3}, though the technical details
(which we omit) are somewhat more complicated, which is related to the
necessary distinction between the cases $s_i=1$ and $s_i=N_i=2$ \cite{fusS3}.
In particular we must now carefully distinguish between Verma and
irreducible modules. Thus in the sequel we will use the previous notations
\vml\ and $\Chil$ only for Verma modules and their \tcha s, while we
write \hill\ for the irreducible quotient of \vml\ and $\chil$ for the
\tcha\ of \hill.
\nextcase
For Verma modules, we find that the combination
  \be  \CHI :=\eE^{-\rho-\Lambda} \Chil \,  \labl{chidef}
(which is independent of $\Lambda$ and corresponds to interpreting
the Verma module as the universal enveloping algebra of $\g_-$)
is odd under the action \erf{w12} of $\Wh$,
  \be  \wh(\CHI)=\epsh(\wh) \,\CHI\, .  \labl{wco}
Here the sign function
  \be \epsh(\wh):=\epsO(\prow^{-1}(\wh))       \labl{newsign}
is the homomorphism $\epsh$ from $\Wh$ to $\zet_2$ that is the pull-back
of the sign function $\epsO$ on $\WO$, rather than the sign function that
$\Wh$ inherits as a subgroup from the sign function $\eps$ of $W$.
\nextcase
For the action of $\Wh$ on the \tcha s of \irmod s we obtain
  \be  \wh(\chil)=\chil  \labl{wcc}
for all $\wh\iN\Wh$, provided that the \hw\ $\Lambda$ is dominant integral;
i.e., the \omchar\ of an \ihwm\ with dominant integral \hw\ is $\Wh$-even.
\internal{
\underline{Proof}:\nextcase
To enter the proof for Verma modules, we choose a specific basis \bminus\
of $\g_-$, including a suitable enumeration of the elements
of \bminus. For $s_i=1$, we choose as the first $N_i$ elements of \bminus\
the step operators $E^{\omD^l i}_- $ for $l=0,1,...\,, N_i-1$
and then the step operators associated to all other negative roots in an
arbitrary ordering. Then according to the \pbw\ theorem the set of all products
$\Enm = \Ene \cdot \Emz$ with $ \Ene:=(E_{}^{-\alpha^{(i)}})^{n_0}_{}\,
(E_{}^{-\alpha^{(\omdi)}})^{n_1}_{}\, \ldots
(E_{}^{-\alpha^{(\omD^{N_i-1}i)}})\raisebox{.7em}{$\scriptstyle n_{N_i-1}$}$
and $\Emz:= (E_{}^{-\beta_1})^{m_1}_{} (E_{}^{-\beta_2})^{m_2}_{} \ldots$,
with $n_i$ and $m_i$ non-negative integers only finitely many of which
are different from zero, forms a basis of $\U(\g_-)$. Commutator
terms that arise when reshuffling the products of generators of $\g_-$
do not contribute to $\CHI$, since $\CHI$ is a trace,
and an element $v^{(\vec n,\vec m)}=\Enm\cdot v_\Lambda^{}$
of the Verma module can contribute to $\CHI$ only
if $n_0 = n_1 =\ldots=n_{N_i-1}^{}=:n$.
The PBW theorem also implies that the contributions to $\CHI$ stemming
from the products $\Ene$ and $\Emz$ factorize, so that
we can investigate their transformation properties under $\whi$ separately.
First, $\whi$ commutes with $\omT$ and maps any $\omT$-orbit of negative roots
to some other orbit of negative roots, i.e.\ only permutes the orbits that
contribute to the second factor, and hence this factor
is invariant under $\whi$. On the other hand, the
contribution $(\CHI)_1$ of operators of the type $\Ene$
to $\CHI$ can be computed explicitly as
$(\CHI)_1 = \eE^{-\rho}/(1-\exp[-\sum_{l=0}^{N_i-1}
\alpha^{(\omD^l i)}])$. Also, evaluating \erf{w12} we find
  \be \whi(\alpha^{(\omD^l i)}) = w_{\omD^l i} (\alpha^{(\omD^l i)})
  = - \alpha^{(\omD^l i)} \, , \labl{wii1}
as well as
  \be \whi(\rho) = \rho - s_i \mbox{$\sum_{l=0}^{N_i-1}$}\alpha^{(\omD^l i)}
  \labl{wref}
for any Weyl vector $\rho\iN\g^\star$, i.e.\ any weight with
$\rho(H^i)=1$ for all $i$. We then obtain
$\whi((\CHI)_1) =\exp[-\rho + \sum_{l=0}^{N_i-1}
\alpha^{(\omD^l i)}]/ (1-\exp[\sum_{l=0}^{N_i-1} \alpha^{(\omD^l i)}])
= - (\CHI)_1$, and hence, combining the two factors, $\whi(\CHI)=-\CHI$.
This holds for all generators $\whi$, and hence we arrive at \erf{wco}.

In the case $s_i=N_i=2$ we choose a different basis \bminus\ of $\g_-$.
As the first three elements of \bminus\ we take the step operators
$E_{}^{-\alpha^{(i)}}$, $E_{}^{-\alpha^{(\omdi)}}$ and $E_{}^{-\alpha^{(i)}
-\alpha^{(\omD i)}}$,
and then again the step operators corresponding to all other negative roots
in an arbitrary ordering. A basis of $\U(\g_-)$ is then given similar as above,
but with the first factor $\Ene$ replaced by
$\Ete:=(E_-^i)^{n_0}_{} (E_-^{\omD i})^{n_1}_{}
(E_{}^{-\alpha^{(i)}-\alpha^{(\omD i)}})^{n'}_{}$.
The contribution to $\CHI$ from operators of the type $\Emz$ again
transforms trivially under $\whi$. Further, in order to have a contribution
$(\CHI)_1$ from operators of the type $\Ete$, we need again $n_0 = n_1
=:n$. The transformation properties $\omega(E_-^i) = E_-^{\omD i}$ and
$\omega( E_-^{\omD i}) = E_-^i$ imply that $\omega( E_{}^
{-\alpha^{(i)}-\alpha^{(\omD i)}}) =- E_{}^{-\alpha^{(i)}-\alpha^{(\omD i)}}$,
which allows us to compute the contribution $(\CHI)_1$ to the \omchar\ as
$(\CHI)_1 =\eE^{-\rho} (1-\eE^{- 2 \alpha^{(i)}-2 \alpha^{(\omD i)}})^{-1}$.
Also, from \erf{w12} we have
  \be \whi(\alpha^{(i)}) = - \alpha^{(\omD i)} \,,
  \quad \whi(\alpha^{(\omD i)}) = - \alpha^{(i)}\,, \quad
  \whi(\alpha^{(i)}+ \alpha^{(\omD i)} ) = -(\alpha^{(i)} +
  \alpha^{(\omD i)}) \, ,\labl{wii2}
while any other positive root is again mapped on a positive root different
from $\alpha^{(i)},\ \alpha^{(\omD i)}$ and $\alpha^{(i)}+ \alpha^{(\omD i)}$.
With \erf{wref} it then follows again that this contribution to $\CHI$ changes
sign under the action of $\whi$, and hence again we obtain \erf{wco}.
\nextcase
For irreducible modules the basic idea is to consider the decomposition
$\calh_\Lambda = \bigoplus \calh_{(L_k)}$ of the irreducible
module \hill\ into irreducible modules of $\g_i$, where
  \be   \g_i := \left\{\bearll \langle E_\pm^{\omdd li}, H^{\omdd li}\mid
  l=0,1,\ldots,N_i-1 \rangle & {\rm for}\ s_i=1\,, \\{}\\[-.8em]
  \langle E_{}^{\pm\alpha^{(i)}\pm\alpha^{(\omD i)}},\, H^i+H^{\omD i}\rangle
  & {\rm for}\ s_i=N_i=2 \,.  \eear \right.\ee
Among the modules $\calh_{(L_k)}$, only those contribute to $\chil$ which are
mapped to themselves by $\tauo$. For $s_i=1$, where $\g_i$
is isomorphic to a direct sum of $N_i$ copies of $A_1$ algebras,
only those basis vectors in these modules contribute to $\chil$ which have
the same weight \wrt all the $A_1$ ideals, $\ell_1=\ell_2=\ldots=\ell_{N_i}$.
These weights are all flipped in sign by $\omT$, and hence $\tauo$ maps a
basis vector $v=v_\ell$ to a vector $v'$ proportional to $v_{-\ell}=(E_-^i
E_-^{\omD i}
\ldots E_-^{\omD^{N_i-1}i})^l\, v$; by the $\om$-twining property \erf{C'}
it then follows that the eigenvalue equation $\tauo(v)= \zeta^k v$ implies
$\tauo(v')=\zeta^k(v')$, i.e.\ $v$ and $v'$ contribute the same phase to
$\chil$. As this is true for all states, it follows that $\whi(\chil)=\chil$,
and as this is true for all values of $i$, we have \erf{wcc}.

In the case $s_i=N_i=2$, where $\g_i$ is isomorphic to $A_1$, there is a
slight complication, because the automorphism $\omega$ acts on $\g_i$ as
$\omega( E_{}^{\pm\alpha^{(i)}\pm\alpha^{(\omD i)}})=
-E_{}^{\pm\alpha^{(i)}\pm\alpha^{(\omD i)}},\;
\omega(H^i+H^{\omD i}) = H^i+H^{\omD i}.$ However,
we now have $v'\propto (E_{}^{-\alpha^{(i)} - \alpha^{(\omD i)}})^{2l} v,$
where the factor of two arises because the $H^i$- and $H^\omdi$-eigenvalues
are added up. Thus only even powers of the step operator
$E^{-\alpha^{(i)}-\alpha^{(\omD i)}}$ occur, so that the additional minus
signs cancel and the reasoning above applies again.
} % endinternal
\nextcase
Next we derive a linear relation between irreducible
and Verma twining characters.  \\[.2em]
If $\mu\leq \lambda$, i.e.\ if $\lambda-\mu= \sum_i n_i\alpha^{(i)}$ with
$n_i\iN\zetpluso$, then
we call $\dpth_\lambda(\mu):= \sum_i n_i$ the depth of $\mu$ \wrt $\lambda$.
By induction on the depth, we can show that
  \be  \CHil = \mbox{{\large$\sum$}$_{\mu\leq\lambda}$}\tilde c_{\lambda \mu}\,
  \chii_\mu^{[\omega]} \, , \labl{551}
with $\CHil$ the \omchar\ of the Verma module with arbitrary
symmetric highest weight
$\lambda$ and $\chii_\mu^{[\omega]}$ the \omchar\ of the irreducible module
with highest weight $\mu$; the numbers $\tilde c_{\lambda \mu}$ are
contained in the cyclotomic field extension ${\dl Q}(\zeta)$ of the rationals
and obey $\tilde c_{\lambda \lambda} =1$.

In order to contribute to the sum \erf{551}, the weights $\mu$ must
obey further requirements in addition to $\mu\le\lambda$.
In particular $\mu$ must be symmetric, and hence we can
restrict the summation to symmetric weights $\mu$ for which $\muO\le\lambdaO$.
Further, the generalized second order Casimir operator
of \g\ has the constant value $C_2(\lambda) = (\lambda+2\rho\Mid\lambda)$
on $V_{\lambda}$; then by \erf{normrel} also
$|\muO+\rhoO|$ has a fixed value, namely $|\lambdaO+\rhoO|$,
for all weights $\mu$ which contribute to \erf{551}.
(Also note that for any  Weyl vector $\rho$ of \g\
the $\gO$-weight $\rhoO\equiv\proj(\rho)$ is a Weyl vector of $\gO$.)
In short, in the decomposition \erf{551} we can restrict $\mu$ to the subset
  \be  \BB(\lambda):= \{\, \mu=\projm(\muO) \mid \muO
  \le\lambdaO, \ |\muO+\rhoO| = |\lambdaO+\rhoO| \,\}  \,.  \ee
Furthermore, for any dominant integral \hw\ $\Lambda$
we can label the elements of $\BB(\Lambda)$ as $\lambda_i$ with $i\iN\natnum$,
in such a way that $\lambdaO_j \le \lambdaO_i$ implies
$i\le j$. Applying the analogue of \erf{551} to all elements
of $\BB(\Lambda)$, it then follows that for all $\lambda_i$ we have
$\Chili = \sum_{\lambda_j\in \BB(\Lambda)} \tilde c_{ij}\, \chilj$
with coefficients $\tilde c_{ij}\iN{\dl Q}(\zeta)$ which satisfy
$\tilde c_{ii} =1$. Further, $\tilde c_{ij}$
can be non-zero only if $\lambdaO_j\le\lambdaO_i$, so that
the matrix $\tilde c=(\tilde c_{ij})$
is upper triangular and can be inverted. Its inverse $c=(c_{ij})$
is upper triangular as well and obeys $c_{ii}=1$, and hence the \omchar\
$\chil$ of the \irmod\ with \hw\ $\Lambda$ can be written as an (infinite)
linear combination
  $\chil = \sum_{\lambda\in \BB(\Lambda)} c^{}_{\lambda} \,\CHil$,
or in terms of the universal Verma twining character $\CHI$,
  \be  \chil = \CHI\cdot\!\! \sum_{\lambda\in \BB(\Lambda)} c^{}_{\lambda} \,
  \eE^{\lambda+\rho}_{}  \,. \labl{lc}
\nextcase
Finally we compare \erf{lc} with the behavior of the
\tcha s under $\Wh$. \\[.2em]
Since $\CHI$ is odd under the action of $\Wh$ while the left hand side of
\erf{lc} is $\Wh$-even, the sum on the \rhs\ must
be $\Wh$-odd, so that $c_{\lambda}=\epsh(\wh)c_{\mu}$
whenever $\wh(\lambda+\rho)=\mu+\rho$ for some element $\wh\iN\Wh$.
Thus for all $\wh\iN\Wh$ we have
$c_\lambda= \epsh(\wh)\, c_{\wh(\lambda+\rho)-\rho}$.
Moreover, with any weight $\lambda$ the weight system of \hill\
already contains the full $\Wh$-orbit of $\lambda$, and hence
we need to know $c_{\lambda}$ only for a single
element of each $\Wh$-orbit, or, equivalently, for a single element of each
$\WO$-orbit of weights of the $\gO$-module \hillO,
say the unique weight in the fundamental Weyl chamber of $\gO$. But the
only such weight $\lambdaO$ which obeys both
$\lambdaO\leq\LambdaO= \proj(\Lambda)$ and $|\lambdaO+\rhoO|=
|\LambdaO+\rhoO|$ is the \hw\ $\LambdaO$ itself; thus only a single
Weyl orbit contributes,
  \be  \chil =  \CHI\cdot \dsum_{\wh\in\Wh} \epsh(\wh)\,\eE^{\wh(\Lambda+\rho)}
                \, . \labl{fast}
Evaluating this identity for the
trivial one-dimensional \irmod\ with highest weight $\Lambda=0$, for which
$ \chii^{[\omega]}_0 = 1$, we have
$\CHI \equiv ( \sum_{\wh\in\Wh} \epsh(\wh)\, \eE^{\wh(\rho)})^{-1}$, so that
  \be \chil =  \sum_{\wh\in\Wh} \epsh(\wh)\, \eE^{\wh(\Lambda+\rho)} \,/
  \sum_{\wh\in\Wh} \epsh(\wh) \, \eE^{\wh(\rho)} \,.   \labl{umf}
When applied to $h\iN\go^{(0)}$, this can be rewritten as
  \be \chil(h) = \frac{ \sum_{\wO\in\WO} \epsO(\wO)\,
  \eE^{(\wO(\LambdaO+\rhoO))(\Pro h)}}
  { \sum_{\wO\in\WO} \epsO(\wO) \, \eE^{(\wO(\rhoO))(\Pro h) }} \,, \ee
and hence
  \be  \chil(h) = \chii_\LambdaO (\Pro h) \,   \ee
by the usual Weyl\hy Kac character formula for integrable highest
weight modules of $\gO$.

This completes the proof of our assertion \erf{1a} for the irreducible \tcha s.
Analogously, the statement for Verma \tcha s
follows by comparing the result for $\CHI$ with the formula for the
Verma module characters of $\gO$ (since $\CHI$ is independent of
$\Lambda$, this latter result holds for arbitrary
highest weights, not just for dominant integral ones).

\sect{Affine algebras}

For the application to \cft, we are interested in the special case where
\g\ is an \uaff\ and where $\om$ is a diagram automorphism that corresponds
to a simple current of the associated \wzwt.\,\,%
\futnote{
These automorphisms can be characterized as the elements of the
unique maximal abelian normal subgroup ${\cal Z}(\g)$ of the
group $\Gamma(\g)$ of diagram automorphisms. This abelian subgroup
is isomorphic to the center of the universal covering Lie group
that has the horizontal subalgebra $\gb\subset\g$ as its \lie.}
In this case we can make a number of additional observations.
\nextcase
In the realization of untwisted \aff s as centrally extended loop algebras
with a derivation, a basis
of \g\ is given by $\HH im$ and $\EE\alphab m$ together with the canonical
central element $K$ and a derivation, where $m\iN\zet$, $i$
takes values in the index set $\Ib$ that corresponds to the \hsa\ \gb\ of \g,
and $\alphab$ is a root of \gb. The full index set is commonly written as
$\I:= \Ib\cup\{0\}=\{0,1,2,...\,,r\}$;\,\,%
\futnote{However, by construction the index set $\ti$ for \gO\ is then
generically {\em not\/} the subset $\{0,1,...\,,\rank\gOb\}$
of $\{0,1,...\,,\rank\gb\}$.}
in this notation, a simple current automorphism must obey $\omd0\ne0$.
%%For the action of the diagram automorphism $\om$ in this basis see
%%\cite{fusS3}.

\internal{
The diagram automorphism $\om$ acts as $\om(K) = K$ and
  \be  \om(\HH in) = \HH{\omd i}n \,, \qquad
  \om(\EE\alphab n) = \etaa\,\EE{\omtb\alphab}{n+\ela} \,. \labl{om}
Here the bar refers to the \hsa\ \gb,
the prefactors $\etaa$ are signs which are +1 for the simple
roots and for all other roots are fixed by the automorphism property of $\om$ ,
$\ela:= (\alphab,\lab{\omdd{-1}0})$
with $\lab i$, $i=\onetor$, the horizontal fundamental weights, and
$\omtb$ is an affine map on the weight space of $\gb$ acting as
  \be  \omtb \lambdab = \kvl\Lambdab_{(\omD 0)} + \sum_{{\scriptstyle
  j=1\atop \scriptstyle j\neq\omD0}}^r \lambdab^{\omD^{-1}j} \Lambdab_{(j)}
  - (\sum_{j=1}^r a_j^\Vee \lambdab^j)\, \Lambdab_{(\omD 0)}\,.  \ee
On the simple \g-roots $\als i$ and the fundamental \g-weights \la i
$\omT$ acts as $\omT(\als i)=\als\omdi$, $\omT(\la i)=\la\omdi$.
} % endinternal
\nextcase
As the generator for the (\onedim) space $\gD$ of derivations one commonly
chooses the element $L_0$ defined by $[L_0, E^i_\pm]=\mp \delta_{i,0}
E^0_\pm$ and $(L_0\Mid L_0)=0$. The automorphism property of $\om$
and the $\om$-invariance of the invariant bilinear form determine
$\omega(L_0)$ uniquely, and there is in fact a unique extension of $\om$
to the semi-direct sum of \g\ and the Virasoro algebra, given by
  \be  \om(L_m)= L_m- ( \lab{\omd0}, H_m) +
  \onehalf\,(\lab{\omd0},\lab{\omd0}) \,\delta_{m,0}\, K \, \labl{omvir}
(and $\om(C)=C$). Using \erf{omvir}, one can show that for any vector of a
\g-module whose weight is symmetric, the action of $\om(L_0)$ coincides with
that of $L_0$; this result is also a consequence of the
Sugawara formula for the Virasoro generators.
\\
The derivation $D$ defined in the general case, which obeys $\omega(D)=D$,
is related to $L_0$ by
  $D = - L_0 + N^{-1} \sum_{l=1}^{N-1}(\Lambdab_{(\omD^l 0)},H)
  - \half N^{-2}\coo\,K$
with $\coo:=\sum_{l,l'=1}^{N-1}(\Lambdab_{(\omD^l
0)},\Lambdab_{(\omD^{l'}0)})$.

\internal{
The relation between the derivations $\LO_0$ (defined analogously as $L_0$)
and $\DO=\Pro(D)$ of $\gO$ is $\DO=-N\LO_0$, which implies
$\Pro(L_0) = N \LO_0 + \Pro(\sum_{l=1}^{N-1} \Lambdab_{(\omD^l0)},H)/N
-\coo\KO/2N^2.$
\nextcase
A sufficient (and, of course, necessary) condition for $\lambda\in\go^\star$
to be a symmetric weight  is $\lambda^i=\lambda^{\omdd li}$
for all $i=\otor$ and all $l$, i.e.\ there is no extra condition on the
$\delta$-component of $\lambda$.
}
\nextcase
Employing the explicit form of the action of $\omT$ on
the simple roots $\alpha^{(i)}$ and fundamental weights \li\ of \g\ and the
action of $\projm$ on the simple roots $\alphaO^{(i)}$ and
fundamental weights $\LambdaO_{(i)}$ of $\gO$, one obtains a series of
nice identities which can be used to simplify formul\ae. For instance,
for the metric of the horizontal part of $\gO^\star$ we can derive the relation
  \be  \overline{\GO}_{ij} = (\Lambdaob_{(i)},\Lambdaob_{(j)})
  = \frac1N\, (\projm\Lambdaob_{(i)},\projm\Lambdaob_{(j)})
  = \frac{N_i N_j}{N^3} \,\Llb \sum_{l,l'=0}^{N-1} \Gb{\omD^l i}{\omD^{l'}j}
      - a_i^\Vee a_j^\Vee \coo \Lrb \,.\labl{Go2}
Thus we have $(\lambdab,\mub) = N\,(\lambdaob,\muob) + \coo\,\kvol\kvom$\,
for $\gO$-weights $\lambdao$ and $\muo$ at levels $\kvol$ and $\kvom$,
respectively.
\nextcase
In particular, the quadratic Casimir eigenvalue of a symmetric
highest \g-weight $\Lambda$ at level $\kv$ can be written as
$\bar C_2(\Lambda)\equiv(\Lambdab,\Lambdab+2\rhob)=N\,(\Lambdaob,
\Lambdaob+2\rhoob)+\coo\,\kvo(\kvo+2\gvo)$.
Further, the dual Coxeter numbers $\gvo\equiv\sumIO i\avO i$
of \gO\ and $\gv\equiv\sumrO i\avi$ of \g\ are related by $N\,\gvo=\gv$, and
for any symmetric \g-weight $\lambda$ of level $\kvl$, the
level of $\lambdaO=\proj(\lambda)$ is $\kvol =\kvl/N$. This implies
a simple relation between the conformal weights
$\Delta_\Lambda\equiv\bar C_2(\Lambda)/((\thetab,\thetab)(\kv+\gv))$
of primary fields of the
$\g$ and $\gO$ \wzwts\ (at levels $\kv$ and $\kvo$, respectively), namely
  \be  \cd_\Lambda = \cdo_{\LambdaO} + \Frac1{2N^2}\,\coo\kV\,
  (1+\Frac{\gV}{\kV+\gV})  = \cdo_{\LambdaO} + \Frac1{24}\, \llb\,
  \Frac{\kV}{\gV}\,(D-\brev D)+c -\brev c \,\lrb    \labl{cd2}
with $D$ and $\brev D$ the dimension of $\gb$ and $\gOb$, respectively.
\nextcase
The order $N$ automorphisms of the \aff s $A_{N-1}\untw$ are generated
by the permutation $\omD$ that `rotates the \dyd\ by one unit'. Thus all
$N$ nodes of the \dyd\ lie on a single $\omD$-orbit, so that
the level of any symmetric integrable highest weight is a multiple of $N$,
and only a single such weight occurs at that level. The prescription
\erf{AO} then formally yields the $1\times1$ `Cartan matrix' $A=0\,$.
\\
The $\om$-invariant subspace $\go^{(0)}$ of $\go$ is \twodim; it is
spanned by the central element $K$ and the derivation $D$.
As $K$ is central, we can therefore write
  \be  \Chil(t,\tau) \equiv \Chil(t K + \tau L_0) = \eE^{2\pi\ii tk_\Lambda}
  \cdot {\rm tr}_{V_\Lambda} \tauo \eE^{2\pi\ii \tau L_0} \, . \labl{73}
Applying the \pbw\ theorem with a suitable basis for $\g_-$, we can
obtain a simple product formula for $\Chil$. By suitably rearranging
factors (within finite products) in this formula, we can then show that
$\ChilTau\equiv\Chil(t=0,q=\eE^{2\pi\ii\tau})$ satisfies the functional
equation $\ChilTau = \ChilnTau$, which implies that $\ChilTau$ is constant.
Evaluating the function at $q=0$, we find that $\ChilTau\equiv 1$, which means
that only the highest weight contributes to the \omchar\ of the Verma module.
As the \hw\ vector is never a null vector, this statement also applies to
the \irmod. Hence
  \be  \chil(t,\tau) = \Chil(t,\tau)
  = \eE^{2\pi\ii tk_\Lambda} \eE^{2\pi \ii \tau \Delta_\Lambda} \, , \labl{cc1}
with $\Delta_\Lambda$ the eigenvalue of $L_0$ on the highest weight vector.
This proves our assertion \erf2.
\nextcase
The characters of $\gO$,
and hence also the \omchar s, have nice modular properties. \\[.2em]
At any fixed non-negative integral value of the level, the set of \ihwm s
with dominant integral \hw s carries a unitary \rep\ of \slz.
This group acts naturally on the modified characters $\tchi_\Lambda\equiv
\eE^{-s_\Lambda \delta} \chii_\Lambda$, where $s_\Lambda$ is the number
  \be s^{}_\Lambda := \Frac1{\raisebox{-.2em}{\small$(\thetab,\thetab)$}}
  \left( \Frac{(\Lambdab+\rhob,\Lambdab+\rhob)}{\kV+\gV} -
  \Frac{(\rhob,\rhob)}{\gV} \right) = \Delta_{\Lambda} - \Frac c{24} \, \ee
(with $c=\kv D/(\kv+\gv)$ the Virasoro central charge), called
the modular anomaly.

In order to inherit the nice modular transformation properties from the
modified characters of \gO, the {\em modified \omchar s} must be defined as
  \be \tchil := \eE^{-\hat s_\Lambda \hat\delta} \chil  \,  \ee
with $\hat\delta=\projm(\deltaO)$ and $\hat s_\Lambda := \sO_{\proj\Lambda}
\equiv \sO_{\LambdaO} = \CO(\LambdaO)/((\brev\ttab,\brev\ttab)(\kvO+\gvO))
- \cO/24$. That is,
the modification relevant to the \omchar\ of \g\ is {\em not\/} the one of the
ordinary character of \g, i.e.\ $\exp(-s_\Lambda \delta)$, but
rather the pull-back of the modification of the $\gO$-character. Thus
  \be  \tchil(h)
  = \exp\llb \sO_{\proj\Lambda} (\projm\,\deltaO)(h)\lrb\, \chiO_{\LambdaO}
  (\Pro h)
  = \exp\llb \sO_{\proj\Lambda}\,\deltaO(\Pro h)\lrb\, \chiO_{\LambdaO}(\Pro h)
  = \tilde {\chiO}_{\LambdaO} (\Pro h) \, . \ee
However, we can show that in fact $\projm(\deltaO)= \delta$, and hence
because of \erf{cd2} the required modification
differs from the modification of the ordinary characters
just by an overall constant which only depends on the level of $\Lambda$:
  \be  s_\Lambda = \cd_\Lambda - \Frac c{24}
  = \cdo_{\LambdaO} - \Frac{\brev c}{24}\, + \Frac1{24}\,\Frac{\kV}{\gV}
  \,(D-\brev D) = \hat s_\Lambda + \Frac1{24}\,\Frac{\kV}{\gV}\,(D-\brev D)
  \,. \labl{cd3}
\nextcasE
For the simple current automorphisms of order two of $\g=C_{2n}^{(1)}$ or
$B_{n+1}\untw$, the \olie\ is the twisted affine \lie\ $\tilde B_n^{(2)}$.\,\,%
 \futnote{We use the notation of \cite{FUch}; in \cite{KAc3}, these
algebras are denoted by $A_{2n}^{(2)}$.}
Our results fit with previous observations at odd level $\kv=2\p+1$ of \g\
\cite{scya6} which suggest a relation of $\gO$ to the $C_n\untw$ WZW theory
at level $\p$. Indeed, we can show that
the modular $S$- and $T$-matrices of $\gO=\tilde B_n^{(2)}$ at level
$\kvo=\kv$ coincide (up to sign factors which are related to certain flips
arising in the application to fixed point resolution) with the
$S$- and $T$-matrices of
$C_n\untw$ at level $\p=(\kv-1)/2$. Also, for even levels our results prove a
conjecture \cite{fuSc2} for the \smat\ which was based on a level-rank duality
of $N=2$ superconformal coset models.

\Sect{Coset theories}{cos}

We can now finally apply what we learned above to the situation outlined in
the beginning, i.e.\ to coset \cfts.

The basic idea of the coset construction (\hsp{-1.4}\cite{goko}, compare also
\cite{baha}) is that, given a reductive \lie\ \gb\ and a reductive subalgebra
$\gbP\emb\gb$, one considers the corresponding embedding $\gP\emb\g$ of
the affinizations \g\ (at some level $\kv\iN\zetpluso$) and \gP\ and
analyzes the difference
  \be  \Lc_m := L_m - \LP_m  \labl{cvir}
of the generators of the two Virasoro algebras which are associated to
the affine \kma s \g\ and \gP\ via the Sugawara formula. The $\Lc_m$ generate
again a Virasoro algebra, with central charge $\cc= c - c'$, the {\em\cvira}.
The basic questions are then the following.
\nextcase
Does this prescription define a consistent \cft\ (which then is referred to
as the {\em coset theory\/} and briefly denoted by `\,$\gb/\gbP$\,')\,?
\nextcase
If so, is that \cft\ unique?
\nextcase
What is the (maximally extended) chiral symmetry algebra \calw\ of the
coset theory?
\nextcase
What is the spectrum of the theory, i.e.\ which \ihwm s (of \calw, or at
least of the \cvira) appear? \\[.2em]
While the coset construction has been proposed already a long time ago,
to a large extent these crucial questions are still open.

Here we will provide an answer to the last question. (This will also
have some impact on all the other questions, but it should be stressed
that in particular the problem of describing the chiral algebra \calw\
is not yet fully solved.) The main concepts needed to arrive at this
answer are branching functions \cite{KAc3}, field identification
\cite{gepn8,levw,mose4}, simple currents \cite{scya,fuge,jf15,scya6},
and fixed point resolution \cite{scya6,scya5,sche3,fuSc2,fusS4},
and it is via the latter concept that \tcha s and \olie s come into the game
\cite{fusS4}.

In more mathematical terms, the question about the spectrum of the coset
theory is the following. Given the \cvira\ associated to a specific pair
\gb\ and \gbP, what are the representation spaces on which it acts?
By construction, the generators \erf{cvir} are defined on the chiral state
space
of the \wzwt\ based on \g, i.e.\ on the direct sum
${\cal H}_\g = \bigoplus_\Lambda \hill$ of all inequivalent
irreducible highest weight modules at level \kv. However, as it turns out,
this is {\em not\/} the state space of the coset theory.
A first crucial observation is that $\gP$ commutes with all the generators
$\Lc_n$. This implies that if we would retain the full state space ${\cal H}_
\g$, then the coset theory would have spin zero fields other than the identity
field, namely all the currents $J^{a'}(z)$ of the subalgebra \gP. To avoid
this disaster, we must require that
these fields act trivially on the state space of the coset theory.
This requirement is implemented by a {\em gauge principle}:
any two states in ${\cal H}_\g$ which differ only by the action of $\gP$
are considered as different descriptions of the same physical
situation, i.e.\ as representing one and the same vector in the
state space of the coset theory.

We conclude that consistency requirements force us to
consider \gP-{\em orbits\/} of vectors in ${\cal H}_\g$ rather than individual
vectors of ${\cal H}_\g$. Concretely, this is implemented by decomposing
the \g-modules $\hl$ into \gP-modules \hLP\ as
 \futnote{We denote quantities referring to the subalgebra \gP\
by the same symbol as the corresponding quantities for \g,
but with a prime added. However, in order to simplify notation, we suppress
`double priming', e.g.\ do not attach a prime to an object that already
has a primed subscript or superscript.}
$  \hl=\bigoplus_\LambdaP \hll\otimes\hLP\,,$
and tentatively considering $\hll$ as irreducible modules of the coset
theory. In terms of characters, this corresponds to regarding
the {\em branching functions} $\bll$ which appear in the decomposition
  \be \chii^\g_\Lambda(\tau) = \sum_{\lP}\bll(\tau)\, \chii^{\gPP}_{\lP}(\tau)
  \labl{cc}
as the characters of the coset theory. However, when doing so
severe problems still remain; they can be attributed to the fact that in
general the redundancy symmetry of the theory is larger than the obvious gauge
symmetry \gP. First, typically there are selection rules, i.e.\ the
branching functions for certain pairs $\ilabel$ vanish. This observation is
in itself not too disturbing, as one just has to make sure to find
all selection rules. (But still, while empirically for most coset theories
all selection rules come from conjugacy class selection rules for the embedding
$\gbP\hookrightarrow\gb$, so far a prescription to
enumerate all selection rules for arbitrary coset theories is not known.)
However, along with each selection rule there also comes a redundancy;
namely, non-vanishing branching functions for distinct pairs $\ilabel$ turn out
to be identical, in particular the putative vacuum module
seems to occur several times. By the same argument which forced us to divide
out the action of the subalgebra \gP, we see that this degeneracy
cannot be interpreted as the multiple appearance of a corresponding primary
field in the spectrum (in particular, it would then
not be possible to obtain the required modular transformation properties).
Rather, the correct interpretation is that a primary field of the coset theory
is not associated to an individual pair $\ilabel$, but rather to an
appropriate equivalence class $\ilabl$ of such pairs. This prescription has
been termed {\em field identification}.

A very convenient description of both conjugacy class
selection rules and field identification is provided by the concept of
{\em identification currents\/} \cite{scya5,scya6}. Namely, first, there is
a subgroup $\Gid$ of the group of integer spin simple currents of the
tensor product theory $\g \oplus (\gP)^*$, the so-called
identification group, such that the conjugacy class selection rules
are given by the condition that the monodromy charges $Q_{(J;J')}$ of any
allowed branching function with respect to all simple currents
$(J;J')$ in the identification group vanish.
Here $Q_{(J;J')}(\ilabel)=Q_J(\Lambda)-Q_{J'}(\Lambda')$, where $Q_J(\Lambda)$
is the combination
$Q_J(\Lambda)=\Delta_\Lambda+\Delta_J-\Delta_{J\star\Lambda}$
of conformal weights.
Second, the equivalence classes in the field identification are precisely
the orbits of the identification group,
  \be  \ilabl = \{\, \jlabel \mid M=J\star\Lambda,\;M'=J'\star\Lambda' \
  {\rm for\ some}\ (J;J')\iN\Gid \,\} \,. \ee

Moreover, provided that \wrt each subgroup of $\Gid$ all orbits
have a common length, the so obtained spectrum corresponds to a consistent
\cft. In particular, by taking one branching function out of each
$\Gid$-orbit and combining them diagonally one arrives at a modular
invariant partition function, and the modular \smat\ is given by
${\cal S}^{}_{\ilabl,\jlabl} =\Sm\Lambda\Mu\,(\SP\LambdaP\MuP)^*$,
where $\ilabel$ and $\jlabel$ are arbitrary representatives of the
orbits $\ilabl$ and $\jlabl$, \resp.

The next degree of generality, and hence difficulty, is reached when orbits
with different sizes appear. (Of course, all sizes are divisors of the size
$N$ of the orbit through $(0;0)$, on which $\Gid$ acts freely,
i.e.\ of the size of $\Gid$. Orbits of non-maximal size are referred to as
{\em fixed points}.) In this situation, taking precisely one
representative out of each orbit $\ilabl$ would not lead to a
modular invariant partition function. Rather, the combination
  \be  Z = \sum_{\scriptstyle [\Lambda;\lP] \atop\scriptstyle Q=0}
  |\Gstab| \cdot |\!\!\! \sum_{(J;J')\in\GId/\Gstab} \!\!\!
  b_{(J;J')\star(\Lambda;\lP)}  \,|^2  \labl{2i}
is modular invariant, where the {\em stabilizer subgroup\/} $\Gstab\subseteq
\Gid$ consists of those elements of $\Gid$ that leave
$(\Lambda;\lP)$ invariant. However, because of $G_{0,0}=\{\id\}$,
$|\Gid|^2$ copies of the vacuum appear in \erf{2i}, so that
one would like to divide \erf{2i} by this factor. But this inevitably leads
to fractional coefficients in the partition function which does not fit with
the interpretation of the partition function as a sum of squares of characters.

A closer look at \erf{2i} suggests
that any fixed point should in fact correspond to $|\Gstab|$ many rather than
to a single primary field. According to such a description, called the
{\em resolution\/} of the fixed point, the labelling of fields
is not by the orbit $\ilabl$ alone, but a fixed point $\ilabel$ splits into
different fields labelled as $(\ilabl,\ch)$ with $\ch$ running from 1
to $|\Gstab|$. As far as \slz-\rep s are concerned, fixed point resolution
amounts to \cite{scya5} decomposing the term in \erf{2i} that corresponds to
a fixed point as the square of a sum of $|\Gstab|$ functions, each of which
differs from the branching function by a certain {\em character modification}
and which is interpreted as the character of an irreducible module.
There is then also a corresponding modification of the \smat\ elements that
involve fixed points.

\internal{
The terms `field identification' and `character modification'
are actually misnomers. Not the fields are to be identified, but rather,
several distinct combinations of highest weights must be identified because
they describe one and the same field of the coset theory. And it is not
the characters which get modified, but rather the branching functions have
to be modified in order that they coincide with
the true characters of the coset \cft.
}

Moreover, in most cases the required
modifications can be interpreted in terms of the characters and \smat\
of another (putative) \cft, called the {\em fixed point theory}.

\sect{Branching spaces and the \facprop}

In the description of fixed point resolution in terms of fixed point \cfts,
part of the information about the coset theory is encoded in the fixed
point theory, so that it is far from obvious how this information
could possibly be obtained from data of $\g$ and $\gP$ alone.
For a large class of coset theories, this puzzle has been resolved in \fssn.
Namely, as described in more detail below, the character modifications can
be expressed as \tbf s associated to \tcha s of $\g$ and $\gP$, and the
fixed point theory corresponds to the coset construction for the
corresponding \olie s $\gO$ and $\gPO$.
This way one solves the fixed point resolution problem not only at the level
of \slz-\rep s, but also directly at the level of the underlying modules
of the chiral algebra.

The class of coset theories we have analyzed are those for which
for each $(J;J')\iN\Gid$ the
diagram automorphisms $\om$ of \g\ and $\omP$ of \gP\ that correspond to the
simple currents $J$ and $J'$, \resp, are related by the {\em\facprop}
  \be  \om\restr{\sgP} =\omP  \,. \labl{omP}
Among the theories satisfying this condition are in particular the
{\em diagonal cosets}, for which \gb\ is the direct sum $\hb\oplus\hb$ of two
copies of a simple \lie\ \hb\ and $\gbP\cong\hb$ is embedded diagonally into
\gb\ (\erf{omP} also holds for generalizations where \hb\ is reductive and
$m$ copies of \hb\ are embedded in $n>m$ copies thereof). Also note
that any two diagram automorphisms $\omega_1$ and $\omega_2$ which correspond
to simple currents commute. We will use automorphisms $\omega_i$
that fulfill the \facprop\ to implement field identification; the fact that
$\omega_1$ and $\omega_2$ commute will imply that
the maps implementing field identification will commute as well.

For diagonal cosets, the embedding of the relevant affine algebras
$\gP \cong \h$ into $\g \cong \h \oplus \h$ is given by
$\EE{\pm\alphab}m= \EE{\pm\alphab}{\eins,m}+\EE{\pm\alphab}{\zwei,m}$,
$\HH im= \HH i{\eins,m}+\HH i{\zwei,m}$ and $K'= K_\eins+K_\zwei$.
Any diagram automorphism respects the triangular decomposition of \g, i.e.\
$\om(\gz)=\gz$ \fortria. (Conversely, any automorphism of an affine \lie\
with this property is induced by a symmetry of the \dyd.)
Also, for diagonal cosets the triangular decompositions of \g\
and \gP\ are compatible in the sense that $\gzP\subseteq\gz$,
and in fact $\gzP=\gz\cap\gP$ \fortria.
When combined with the $\omega$-twining property \erf{C'} of $\tauo$,
the \facprop\ \erf{omP} then implies that also $\omP(\gzP)\subseteq\gzP$.

Consider now a unitary highest weight module $V$ of \g\
on which the central elements $K_\eins$ and $K_\zwei$ act as non-negative
integer multiples $k_\eins,\, k_\zwei$ of the identity.
We can view $V$ as a $\gP$-module, on which the central element $K'$
acts as the multiple $k'=k_\eins + k_\zwei$ of the identity.
The coset Virasoro algebra then acts on $V$ via
the generators $\Lc_m=L_{(1),m}+L_{(2),m}
-L'_m$, where the individual terms are obtained via the Sugawara construction
for \h\ at levels $\KV_\eins=2k_\eins/(\ttab\Mid\ttab),\,
\KV_\zwei=2k_\zwei/(\ttab\Mid\ttab)$ and $\KV_\eins+\KV_\zwei$, \resp.
The identification group for diagonal cosets consists of all simple currents
of the form $(J;J')=(\KV_\eins\Lambda,\KV_\zwei\Lambda;(\KV_\eins+\KV_\zwei)
\Lambda)$ with $\Lambda$ a so-called cominimal \cite{fuge} fundamental
\h-weight, and hence is isomorphic to the group of those symmetries
of the Dynkin diagram of \h\ which do not leave the zeroth node fixed.

Moreover, because of $\omP(\gpP)\subseteq\gpP$ we have
$\gpP\Cdot(\tauo(v))=0$ whenever $\gpP\Cdot v=0$; thus $\tauo$ maps \hwv s
\wrt \gP\ in \hl\ to \hwv s \wrt \gP\ in ${\calH}_{\omtLa}$.
Similarly, if $v_1$ and $v_2$ lie in one and the same irreducible \gP-module,
then so do $\tauo(v_1)$ and $\tauo(v_2)$.
It follows that on any irreducible \g-module \hl\ the map $\tauo$ factorizes
into a mapping $\tauc$, of the same order $N$, which maps any irreducible
\gP-submodule of \hl\ to some \gP-submodule of \hlo\
(i.e., a map for which, roughly speaking, only the positions of the \gP-modules
$\hlP$ and $\hlPt$ in the respective \g-modules matter) and the map $\tauoP$
defined via the automorphism $\omP$ of \gP.
Hence, decomposing \hl\ into \gP-modules as
  \be  \hl=\bigoplus_\el \hlP  \, \labl{hh}
(thus the \gP-modules in \hl\ are labelled by $\el=1,2,...$, and as a
redundant information also the highest \gP-weight $\lP\equiv\lP(\el)$ is
displayed), the analogous decomposition of $\hlo$ reads
  \be  \hlo\equiv\tauo(\hl)=\bigoplus_\el \tauo(\hlP)
  =\bigoplus_\el \hlPt  \,. \labl{hho}

In terms of the branching spaces
  \be  \hll= {\rm span}_\compleX\{\vlaP\mid\lambdaP(\el)=\LambdaP\}
  \labl{hll}
(\vlaP\ denotes the \hwv\ of $\hlP$), the decomposition \erf{hh} reads
$\hl=\bigoplus_\el \hlP \cong \bigoplus_\LambdaP \hll\otimes\hLP$,
where the summation is over all integrable highest weights of \gP\ at the
relevant level.
The branching functions $\bll$ appearing in \erf{cc} are the traces
  \be \bll(\tau) = \trhll \qc \,  \ee
($q=\exp(2\pi\ii\tau)$) over the branching spaces.
{}From the \facprop\ it follows that
  \be  b_{\oMT\Lambda,\oMT'\Lambda'} = b_{\Lambda,\Lambda'} \,,\labl{blb}
in agreement with the fact that $(\om;\omP)$ corresponds to an
identification current $(J;J')$.

Further, it also follows from \erf{omP} that the \cvira\ is $\om$-invariant,
  \be  \om(\Lc_m) =\om(L_m)-\omP(\LP_m) = \Lc_m \,.  \labl{omcvir}
Together with the $\omega$-twining property \erf{C'} of $\tauc$,
\erf{omcvir} implies that $\tauc\!:\;\hll\to\hlolo$
intertwines the action of the \cvira, i.e.\ we have
  \be  [\tauc,\Lc_m]=0 \,   \labl{tclc2}
for all $m\iN\zet$. In particular, $\hll$ and $\hlolo$ carry isomorphic \rep s
of the coset Virasoro algebra. Now in order to keep the interpretation of
\gP\ as a redundancy symmetry, not only the \cvira, but also
the full (maximally extended) chiral algebra \calw\ of the coset theory
must be well-defined on \gP-orbits, i.e.\ \calw\ must commute with \gP.
(As a consequence, $\hll$ and $\hlolo$ also carry isomorphic \rep s
of \calw, as is already suggested (but not implied) by \erf{blb}.) Further,
compatibility with field identification requires that the intertwining property
of $\tauc$ must be valid not only for the \cvira, but for \calw\ as well.

\sect{Fixed points}

If $\ilabel$ is a fixed point of an identification current, i.e.\ if
$\omt(\Lambda)=\Lambda$ and $\omtP(\LambdaP)=\LambdaP$,
the map $\tauc$ is an {\em endo}morphism, so that we can
define the {\em\tbf} $\bllo$ as the trace
  \be  \bllo(\tau) := \trhll\, \tauc \, \qc   \,. \labl{bllo}
(For diagonal cosets a fixed point is a combination
$(\Lambda\eins\,,\Lambda\zwei\,;\LambdaP)$ of three weights of \h\ that are
fixed \wrt the same simple current automorphism $\omh$ of \h).
Using the fact that $[\Lc_m,L'_n]=0$, so that $q^{L_m}=
q^{L'_m}\cdot q^{\Lcs_m}$, and that $[\tauc,\Lc_m]=0$, we obtain
  \be  \chil(\tau)= \sum_{\LambdaP} \bllo(\tau)\,\chilP(\tau) \,, \labl{cbc}
i.e.\ $\bllo(\tau)$ arises in the decomposition of the \ocha s of \g\
\wrtt automorphism $\tauo$ into \ocha s of \gP\ \wrt $\tauoP$
(this justifies the name `\tbf'). In particular, analogously as for
ordinary branching functions, the \tbf s can in principle be computed from the
\tcha s $\chil$ and $\chilP$. While doing this directly would be cumbersome,
one can employ the equivalence between \tcha s and the
characters of the \olie s, so that \bllo\ can be computed as the ordinary
branching function
    \be  \bllo = \brev b^{}_{\brev\Lambda,\brev\Lambda'} \labl{brevb}
of the coset construction $\gO/\gO{}'$ of \olie s.  (Above  only the case
when \gb\ and $\gOb$ are simple has been described explicitly; in the
general case, one deals with a product of \omchar s of the various
summands of \g\ \resp\ \gP.)

As a consequence, the \tbf s possess nice modular transformation properties.
In particular, under $\tau\mapsto-1/\tau$ they transform with
  \be  \Sob\Lambda\LambdaP\Mu\MuP= \So\Lambda\Mu\,(\SoP\LambdaP\MuP)^*
  = \SoO\LambdaO\MuO\,(\SoOP\lambdaOP\muOP)^* \,, \labl{Sob}
where $\So\lambda\mu$ \resp\ $\SoP\lambdaP\muP$ are the
corresponding $S$-matrices for the \omchar s of \g\ \resp\ \gP,
which in turn are nothing but the $S$-matrices of the respective \olie s.
Also, up to an over-all phase the matrix ${\cal T}^{[\omega]}$ which describes
the behavior under
$\tau\mapsto\tau+1$ coincides with the restriction of the ordinary (untwined)
$T$-matrix $T^{\g}\otimes (T^{\g'})^*$ to the fixed points.

It is also not difficult to see that if $J$ is a simple current of \g,
then $\brev J$ is a simple current of \gO\ (albeit possibly the trivial
one, i.e.\ the identity field). For coset theories, an
analogous statement holds for identification currents, and hence each
selection rule of $\g/\gP$ leads to a selection rule of $\gO/\gO{}'$.
As a consequence, the matrix element $\Sob\Lambda\LambdaP\Mu\MuP$ does not
depend on the choice of representatives of the orbits of $\ilabel$ and
$\jlabel$. This implies e.g.\ that, analogously to \erf{blb},
  \be  b^{[\om_1]}_{\om_2^\star\Lambda,\om_2^\star{}'\Lambda'}
  = b^{[\om_1]}_{\Lambda,\Lambda'} \,.\labl{bob}

Now recall from the discussion after \erf{tclc2} that field identification
requires that any element of the coset chiral algebra \calw\ commutes with
$\tauc$. For fixed points this implies that already
the eigenspaces \wrt $\tauc$ of the branching spaces form modules over
the \cvira, and in fact over the coset chiral algebra $\calw$.
This way the implementation of field identification by the
maps $\tauc$ enforces fixed point resolution: one has to split the branching
spaces \hll\ into the eigenspaces \wrt the action of the stabilizer group
$\Gilabel$. Thus the label $\ch$ of the resolved fixed points $(\ilabl,\ch)$
corresponds to the group characters $\ch\!:\, \Gilabel\to\complex$
of $\Gilabel$; each $\Gilabel$-character $\ch$ gives rise to an eigenspace
\hllch\ of \hll\ on which $\tauc$ has eigenvalue $\ch(\om)$.
Thus, denoting the (Virasoro specialized) character of the eigenspace
\hllch\ by \bllch, we have
  \be  \bll(\tau) = \sum_{\ch\in\Gstabs} \bllch(\tau) \,, \labl{bb}
where $\Gilabels$ is the group of characters of $\Gilabel$. Since
the coefficients in the expansion of \bllch\ in powers of $q$
count the number of eigenstates at the relevant grade, they are
manifestly non-negative integers, and hence the functions \bllch\
can be interpreted as the characters of the resolved fixed points.
Further, we have $\bllo = \sumgilabels\!\! \ch(\om)\,\bllch$, which by
the orthogonality of the group characters can be inverted to
  $\bllch = \sumgilabel \!\!\ch^*(\om)\, \bllo /|\Gilabel|$.
Of course, one still has to
implement the field identification \wrtt rest of \Gid, which however,
owing to \erf{bob}, just amounts to considering orbits of identical
\tbf s. Thus, employing also the convention to set $\ch(\omega) = 0$
for $\ch\iN\Gilabels$ whenever $\omega\not\in \Gilabel$, the characters of the
resolved fixed points are
  \be  \bearll \chillch \!\!&= \Frac1{|\GId|} \dsum_{\omega\in\GId}\!
  b^{(\ch)}_{\omt\Lambda,\omtP\Lambda'}
  = \Frac1{|\GId|\cdot |\Gilabel| } \dsum_{\ome,\omz\in\GId} \ch^*(\ome) \,
  b^{[\ome]}_{\omzD\Lambda,\omzDP\Lambda'} \\{}\\[-.5em]
  & = \Frac1{|\Gilabel|} \dsum_{\om\in\GId} \ch^*(\om) \,
  b^{[\om]}_{\Lambda,\Lambda'} = \bllch \, . \eear \labl{crF}
In other words, while in the presence of fixed points
the true characters of a coset conformal field theory
are not the branching functions, they can be obtained as
linear combinations of branching functions of the original coset theory
and those for the embedding of the associated \olie s.
Also note that when combining the two chiral halves in the diagonal modular
invariant, we obtain less combinations of states than one would naively
expect. This is the projection alluded to in the first section.

For the \smat\ which describes the
behavior of these characters under the modular transformation
$\tau\mapsto-1/\tau$, we find (by first expanding $b^{[\om]}_{\Lambda,\Lambda'}
(-1/\tau)$ in terms of the $b^{[\om]}_{\Mu,\Mu'}(\tau)$ and then re-expressing
the $b^{[\om]}_{\Mu,\Mu'}$ through the $\bllch$)
  \be  \Scc\Lambda\LambdaP M{M'} =
  \frac{|\GId|}{|\Gilabel|\cdot |\Gjlabel| } \sum_{\om\in\GId}
  \ch^*(\om)\, \CS^{[\om]}_{\ilabl,\jlabl} \, \cht(\om) \,. \labl{Scc}
Note that $\ch\iN\Gilabels$ and $\cht\iN\Gjlabels$, so that
the summation is effectively only over the subgroup
$\Gilabel\cap\Gjlabel\subseteq\Gid$ of the identification group. Thus
if either $\ilabel$ or $\jlabel$ is not a fixed point, the sum in \erf{Scc}
contains only a single term, namely the one with $\om=\id$, and the \smat\ is
just a multiple of the ordinary \smat\ element for the corresponding branching
functions.
Similarly, the transformation under $\tau\mapsto\tau+1$ is implemented by
  \be \CT_{(\ilabl,\ch),(\jlabl,\cht)} = \delta_{\ilabl,\jlabl}\,
  \delta_{\ch,\cht}\, T_{\ilabel} \, . \labl{tdef}
One can check that the matrix \erf{Scc} is symmetric and unitary,
that its square gives a permutation of order two on the primary fields of the
coset \cft, and that $\CS$ and $\CT$ generate a unitary \rep\ of \slz.

As an illustration, consider the case $\Gid\cong\zet_N$ with $N$ prime. Then
$\Gid$ has $N$ characters $\ch=\ch_k$, $k=\otoNm$, acting on $n\iN\zet_N$
as $\ch_k(n) = \zeta^{kn}$ with $\zeta=\exp(2\pi\ii/N)$. In this case for
$\om\iN\Gid$ ($\om\ne\id$) the matrix $\SNb$ is the same for all $n\ne0$,
  $\Snb\Lambda\LambdaP\Mu\MuP =: \Sf\Lambda\LambdaP\Mu\MuP$, while
  $\Sno\Lambda\LambdaP\Mu\MuP = \Soo\Lambda\LambdaP\Mu\MuP
  \equiv \Sm\Lambda\Mu(\SP\LambdaP\MuP)^*$.
In the interesting case when both combinations of weights
are fixed points, we then obtain
  \be  \Skl\Lambda\LambdaP\Mu\MuP =
  \Frac1N\, \Soo\Lambda\LambdaP\Mu\MuP + \llb \delta^{}_{k,l}-\Frac1N \lrb\,
  \Sf\Lambda\LambdaP\Mu\MuP  \,, \ee
which is precisely the expression for the \smat\ that had been conjectured
in \cite{scya5}.

\sect{Further applications and problems}

We conclude with some remarks on related issues and open questions.

\nextcomm
On the purely mathematical side, one may wonder whether the notion of
orbit \lie\ can still be given a meaning for diagram automorphisms which
do not satisfy the \ocond\ \erf{ocond}, maybe as a Borcherds algebra
or as a generalized \kma\ in the sense of \cite[\S 11.13]{KAc3}.
\\
Also, we believe that the dual \lie\ of \gO, i.e.\ the \lie\ obtained from
\gO\ by reversing the direction of the arrows of the \dyd, is precisely
the subalgebra of the dual algebra of \g\ that stays fixed under
the automorphism $\tilde\om$ of the dual \lie\ that corresponds to $\om$.\,\,%
\futnote{In the case of affine algebras, these dual \olie s
have been used in the reduction of Toda field equations \cite{oltu,bcds2}.}
A proof of this connection would certainly be welcome.

\nextcomm
In the case of \uaff s, the characters of the orbit \lie\ $\gO$, and hence
the \omchar s as well, carry a unitary \rep\ of $SL(2,\zet)$ whenever
the diagram automorphism corresponds to a simple current, i.e.\ belongs to
the maximal abelian normal subgroup ${\cal Z}(\g)$ of the
group of diagram automorphisms. This should be related to the fact
\cite{fugv2} that the automorphisms in ${\cal Z}(\g)$ are precisely the
{\em localizable} diagram automorphisms.

\nextcomm
The \olie s of \uaff s correspond precisely to the fixed point theories
introduced in \cite{scya5}, except for the cases of the order two
automorphisms of $\g=C_{2n}^{(1)}$ or $B_{n+1}\untw$, where the \olie\ is
the twisted algebra $\tilde B_n^{(2)}$. In the latter cases so far
the \olie\ could not be interpreted in terms of a genuine
\cft. From the results of \cite{scya6} one also knows that such a \cft\
would have to be non-unitary.

\nextcomm
In contrast to (generalized) diagonal cosets, for general coset theories
the \facprop\ \erf{omP} does not hold. Thus in the general case one
cannot restrict oneself to diagram automorphisms $\omega,\, \omega'$ of
\g\ and \gP. Rather one has to allow for more general outer automorphisms
of \g\ and \gP, which are the product of a diagram automorphism
and some compensating inner automorphism corresponding to an element of the
Weyl group.  However, we expect that the generalization of our ideas to
arbitrary cosets will not involve any further new \rep\ theoretic concepts,
but rather only
additional technical complications coming from the fact that $\omP$ and
$\om\restR{\sgP}$ differ by an inner automorphism. In particular, we expect
that the formula \erf{Scc} for the \smat\ of the coset theory will still be
valid.

\nextcomm
In all cases that we have checked explicitly, the matrix $\CS$
defined by \erf{Scc} gives rise, via the Verlinde formula,
to non-negative integral fusion coefficients. But a general proof of this
property is still lacking.

\nextcomm
It is worth stressing that in the description of the coset construction,
there are still open conceptual problems which are
largely independent of the issues of field identification and fixed point
resolution. Most importantly, it has not yet been proven
that there exists a suitable closure
$\ugtw$ of the universal enveloping algebra $\ug$ such that the
algebra \calw\ defined with the help of $\ugtw$
fulfills all requirements for a chiral algebra of a \cft.
A candidate for $\ugtw$ is the topological completion of $\ug$ that is
described in \cite{frzh}.
\\
In order to prove that \calw\ is indeed the -- maximally extended --
chiral algebra of the coset theory, one must in particular show that
(in the absence of field identification fixed points) the branching
spaces are {\em irreducible\/} modules of \calw, and that each isomorphism
class of irreducible modules appears precisely once. (Note that
as modules over the \cvira, the branching spaces are highly reducible
as soon as $\cc\geq1$.) In the presence of fixed points, one has to
prove the same statement for the eigenspaces of the maps $\tauc$
(in principle a further splitting of these eigenspaces could be necessary,
but we believe that they
are in fact already irreducible modules over $\calw$).
\\
Note that according to our construction the candidate
coset chiral algebra \calw\ is {\em not\/} obtained as the commutant \calwc\ of
$\gP$ in \ugtw, but rather as the subalgebra of \calwc\ that is fixed by all
automorphisms yielding field identifications; this reflects the fact that not
only the action of the currents of \gP, but also the action of the maps
$\tauc$ connects different representatives of one and the same physical state.
More specifically, for a field identification automorphism $\om$ the \facprop\
implies that $\omega(\calwc)\subseteq \calwc$, so that $\omega$ induces an
order-$N$ automorphism of \calwc\ and we can decompose \calwc\ into its
eigenspaces with respect to $\om$. The eigenspace $\calwoo$ that is left
invariant under all these automorphisms is a subalgebra of \calwc\ and
contains all elements of $\calw$ that are intertwined by $\tauc$
(in particular, according to \erf{omvir}, the coset Virasoro algebra).
The fact that field identification is implemented by the maps $\tauc$ then
implies that $\calw = \calwoo$.

\nextcomm
The interpretation of \gP\ as a gauge symmetry is a basic ingredient in the
geometric approach to coset theories via WZW sigma models.
There are indications \cite{hori} that in this formulation
fixed point resolution can be implemented
by considering also topologically non-trivial ${\sf G}'$-principal bundles,
where ${\sf G}'$ is the relevant \Lie\ whose \lie\ is \gP. Our results
therefore suggest interesting new
relations between these bundles in different topological sectors.

\nextcomm
Twining characters and \olie s are natural structures associated to diagram
automorphisms of symmetrizable \kma s. We believe that
they should be useful in various different situations where such
automorphisms play a
r\^ole. As for the case of \cft, the following applications come to mind:
the explanation of the fixed point resolution that is present \cite{scya6}
in integer spin simple current invariants, the influence of fixed points
on the relation \cite{pezu} between \cft\ and graphs, and the computation
of the generalized $S$-matrices $S(J)$ that were defined in \cite{mose3}.

\nextcomm
Finally let us remark on the issue of fixed point resolution in more
general situations.
The fact that we could attribute the presence of different orbit lengths
to the action of a {\em finite\/} group of redundancy symmetries
certainly constitutes a significant simplification. In particular, it
enabled us to describe fixed point resolution in terms of the
eigenspaces of these redundancy symmetries. This will typically no longer
be possible in the general case.
\\
Another major advantage of the theories considered here is
that it is easy to see that of a naive implementation of the redundancies
according to \erf{fig1} leads to inconsistencies. In other contexts such
inconsistencies, if present, will be much more difficult to detect.

%\newpage
\vskip3.3em

 \newcommand\wb       {\,\linebreak[0]} \def\wB {$\,$\wb}
 \newcommand\Bi[1]    {\bibitem{#1}}
 \newcommand\Erra[3]  {\,[{\em ibid.}\ {#1} ({#2}) {#3}, {\em Erratum}]}
 \newcommand\BOOK[4]  {{\em #1\/} ({#2}, {#3} {#4})}
 \newcommand\vypf[5]  {{#1} [FS{#2}] ({#3}) {#4}}
 \renewcommand\J[5]   {{#1} {#2} ({#3}) {#4} }
 \renewcommand\J[5]   {\ {\sl #5}, {#1} {#2} ({#3}) {#4} }
 \newcommand\Prep[2]  {{\sl #2}, preprint {#1}}
 \def\jf    {J.\ Fuchs}
 \def\anop  {Ann.\wb Phys.}
 \def\foph  {Fortschr.\wb Phys.}
 \def\hepa  {Helv.\wb Phys.\wB Acta}
 \def\ijmp  {Int.\wb J.\wb Mod.\wb Phys.\ A}
 \def\jopa  {J.\wb Phys.\ A}
 \def\npbF  {Nucl.\wb Phys.\ B\vypf}
 \def\npbp  {Nucl.\wb Phys.\ B (Proc.\wb Suppl.)}
 \def\nuci  {Nuovo\wB Cim.}
 \def\nupb  {Nucl.\wb Phys.\ B}
 \def\phlb  {Phys.\wb Lett.\ B}
 \def\comp  {Com\-mun.\wb Math.\wb Phys.}
 \def\lemp  {Lett.\wb Math.\wb Phys.}
 \def\phrd  {Phys.\wb Rev.\ D}
 \def\mpla  {Mod.\wb Phys.\wb Lett.\ A}
 \def\duke  {Duke\wB Math.\wb J.}

 \def\A       {Algebra}
 \def\alg     {algebra}
 \def\Be     {{Berlin}}
 \def\BIR    {{Birk\-h\"au\-ser}}
 \def\Ca     {{Cambridge}}
 \def\con     {conformal\ }
 \def\CUP    {{Cambridge University Press}}
 \def\furu    {fusion rule}
 \def\GB     {{Gordon and Breach}}
 \def\fts     {field theories}
 \def\ide     {identification}
 \def\Infdim  {Infinite-dimensional}
 \def\NY     {{New York}}
 \def\nn      {$N=2$ }
 \def\oa      {operator algebra}
 \def\Q       {Quantum\ }
 \def\q       {quantum\ }
 \def\qg      {quantum group}
 \def\Rep     {Representation}
 \def\SV     {{Sprin\-ger Verlag}}
 \def\syms    {sym\-me\-tries}
 \def\va      {Vira\-soro algebra}
 \def\wzw     {WZW\ }

\end{document}